\input amstex
\input epsf
\documentstyle{amsppt}\nologo\footline={}

\font\oitossi=cmssqi8
\font\oitoss=cmssq8
\font\bff=cmbx12 at 14pt

\def\T{\mathop{\text{\rm T}}}
\def\Lin{\mathop{\text{\rm Lin}}}
\def\GL{\mathop{\text{\rm GL}}}
\def\Gr{\mathop{\text{\rm Gr}}}

\def\opensub{\mathop{\subset\!\!{\circ}}}
\def\opensup{\mathop{{\circ}\!\!\supset}}
\def\U{\mathop{\text{\rm U}}}
\def\ta{\mathop{\text{\rm ta}}}
\def\area{\mathop{\text{\rm area}}}
\def\Arg{\mathop{\text{\rm Arg}}}
\def\Re{\mathop{\text{\rm Re}}}

\def\Z{\mathop{\text{\rm Z}}}
\def\I{\mathop{\text{\rm I}}}
\def\B{\mathop{\text{\rm B}}}
\def\S{\mathop{\text{\rm S}}}
\def\E{\mathop{\text{\rm E}}}
\def\arccosh{\mathop{\text{\rm arccosh}}}

\input lpic-mod
\input epsf

\hsize450pt\topmatter\title\endtitle\author\endauthor\address
\endaddress\email\endemail
%\subjclass\endsubjclass
\abstract These lecture notes are based on [arXiv: math/0702714,
0907.4469, 0907.4470]. We introduce and study basic aspects of
non-Euclidean geometries from a coordinate-free viewpoint.
\endabstract\endtopmatter\document

\count0=0

\newpage

\centerline{\bff Basic coordinate-free non-Euclidean geometry}

\vskip15pt

\centerline{\smc Sasha Anan$'$in and Carlos H.~Grossi}

\vskip30pt

\leftskip55pt\rightskip55pt

\noindent
{{\bf1.~Darwin and geometry}\hfill3}

{1.1.~Risks of travelling around the world\hfill3}

{1.2.~Question\hfill4}

{1.3.~Cogito, ergo sum\hfill4}

{1.4.~Guide to the reader\hfill4}

\smallskip

\noindent
{{\bf2.~Projective spaces and their relatives}\hfill5}

{2.1.~Problem\hfill5}

{2.6.~Projective space\hfill7}

{2.7.~Sphere and stereographic projection\hfill8}

{2.9.~Riemann sphere\hfill9}

{2.12A.~Grassmannians\hfill9}

\smallskip

\noindent
{{\bf3.~Smooth spaces and smooth functions}\hfill9}

{3.1.~Introductory remarks\hfill10}

{3.2.~Sheaves of functions\hfill10}

{3.3.~Basic example\hfill10}

{3.4.~Smooth maps and induced structures\hfill11}

{3.5.~Product and fibre bundle\hfill12}

{\quad3.5.1.~Product\hfill12}

{\quad3.5.5.~Fibre product\hfill12}

{3.6.~Tangent bundle\hfill13}

{\quad3.6.1.~Tangent vectors\hfill13}

{\quad3.6.3.~Differential\hfill13}

{\quad3.6.5.~Tangent bundle of a subspace\hfill14}

{\quad3.6.8.~Equations\hfill14}

{\quad3.6.10.~Taylor sheaves\hfill15}

{\quad3.6.11.~Prevarieties\hfill16}

{3.7.~$C^\infty$-manifolds\hfill18}

{\quad3.7.2.~Families and bundles\hfill18}

{\quad3.7.8.~Tangent vector to a curve\hfill18}

{3.8A.~Final remarks\hfill19}

\smallskip

\noindent
{{\bf4.~Elementary geometry}\hfill19}

{4.1.~Some notation\hfill19}

{4.2.~Tangent space\hfill20}

{4.3.~Metric\hfill20}

{\quad4.3.2.~Length and angle\hfill21}

{4.4.~Examples\hfill21}

{4.5.~Geodesics and tance\hfill22}

{\quad4.5.4.~Spherical geodesics\hfill22}

{\quad4.5.6.~Hyperbolic geodesics\hfill23}

{\quad4.5.7.~Triangle inequality\hfill23}

{\quad4.5.8.~Duality\hfill24}

{4.6.~Space of circles\hfill24}

{4.7.~Complex hyperbolic zoo\hfill25}

{4.8.~Finite configurations\hfill25}

{\quad4.8.1.~Stollen Carlos' lemma\hfill26}

{4.9.~There is no sin south of the equator\hfill27}

{4.10A.~Geometry on the absolute\hfill28}

{4.11.~A bit of history\hfill29}

{\quad4.11.1.~References\hfill30}

\smallskip

\noindent
{{\bf5.~Riemann surfaces}\hfill30}

{5.1.~Regular covering and fundamental group\hfill30}

{5.2.~Discrete groups and Poincar\'e polygonal theorem\hfill30}

{5.3.~Teichm\"uller space\hfill30}

\smallskip

\noindent
{{\bf6.~Appendix: Largo al factotum della citta}\hfill31}

{6.15.~Gram-Schmidt orthogonalization\hfill32}

{6.20.~Sylvester's law of inertia\hfill33}

{6.21.~Sylvester criterion\hfill34}

\smallskip

\noindent
{{\bf7.~Appendix: Basic algebra and topology}\hfill35}

\smallskip

\noindent
{{\bf8.~Appendix: Classification of compact surfaces}\hfill35}

\smallskip

\noindent
{{\bf9A.~Appendix: Riemannian geometry}\hfill35}

\smallskip

\noindent
{{\bf10A.~Appendix: Hyperelliptic surfaces and Goldman's
theorem}\hfill35}

\smallskip

\noindent
{{\bf Hints}\hfill36}

{\quad1.2.\hfill36}

{\quad2.3.\hfill36}

{\quad2.5.\hfill36}

{\quad2.10.\hfill36}

{\quad2.11.\hfill36}

{\quad3.3.2.\hfill36}

{\quad3.3.3.\hfill36}

{\quad3.3.4.\hfill36}

{\quad3.6.6.\hfill36}

{\quad4.3.3.\hfill36}

{\quad4.4.1.\hfill36}

{\quad4.5.2.\hfill36}

{\quad4.5.9.\hfill36}

{\quad4.8.2.\hfill37}

{\quad6.5.\hfill37}

{\quad6.6.\hfill37}

{\quad6.7.\hfill37}

{\quad6.9.\hfill37}

{\quad6.10.\hfill37}

{\quad6.23.\hfill37}

{\quad6.24.\hfill37}

{\quad6.25.\hfill37}

\leftskip0pt\rightskip0pt

\newpage

\leftheadtext{Basic coordinate-free non-Euclidean geometry}

\centerline{\bff Basic coordinate-free non-Euclidean geometry}

\vskip15pt

\centerline{\smc Sasha Anan$'$in and Carlos H.~Grossi}

\vskip36pt

\rightline{\oitossi The introduction of numbers as coordinates $\dots$
is an act of violence $\dots$}

\medskip

\rightline{{\oitoss --- HERMANN WEYL,} {\oitossi Philosophy of
Mathematics and Natural Science}}

\bigskip

\rightline{\oitossi`You see, the earth takes twenty-four hours to turn
round on its axis---'}

\rightline{\oitossi `Talking of axes,' said the Duchess, `chop off her
head!'}

\medskip

\rightline{{\oitoss --- LEWIS CARROLL,} {\oitossi Alice's Adventures in
Wonderland}}

\vskip60pt

\bigskip

\centerline{\bf1.~Darwin and geometry}

\medskip

The subject of this course is deeply related to three great geometers:
Riemann, Klein, and Poincar\'e. With the study of the spherical and
hyperbolic plane geometries,\footnote{It turns out that Euclidean
geometry is degenerate and somehow separates the other two.}
we try to achieve the modest purpose of illustrating some contributions
of these geometers. On the way, we encounter a few tools discovered in
recent years.

Human geometric intuition is considerably stronger than the algebraic
one for an obvious reason: since the Stone Age, we have lots of
experience in space moving and very few in counting.\footnote{It was
common for a Neolithic man to keep hunting new wives, not remembering
how many were already in his cave. The uprise of monogamy as a solution
to this problem shows the difficulties with arithmetics at those ages.}

The plane geometries have a strong influence in modern geometry which
can be partially explained on a biological basis. Birds are certainly
excellent geometers: it suffices to see how they express their
(probable) happiness with sophisticated pirouettes in three dimensions
when the rain is over. Serpents ought to be good topologists. (Soon, we
will study a little bit of topology.) Unfortunately, the experience of
human beings is nearly two-dimensional, at most $2.5$-dimensional.
Believers in Darwin's theory might have inherited the three-dimensional
experience from the apes but we doubt that this theory\footnote{As a
believer in evolutionism and its new branches, Carlos does not share
this view. Nevertheless, he is surprised with having difficulties in
combinatorics in spite of his microbiological past.}
actually works: we have never seen an ape turning into a man and are
tired of seeing how it goes the other way around. Hence, in order to be
well armed for our future, it is essential to study geometry. (See as
an illustration the webpage http://www.ihes.fr/$\sim$gromov of one of
the greatest geometers of our days, Misha Gromov.)

\medskip

{\bf1.1.~Risks of travelling around the world.} In three dimensions, a
human being usually possesses two legs. This seems to be enough,
although we would fall less frequently if they were three. Therefore,
in a two-dimensional world, one leg should suffice. Say, the right one.

In the best of all possible two-dimensional worlds, Candid, a son of a
cheerful mother, decided to pursuit a most challenging adventure, to
travel around the world. Fearing the dangers of the voyage, the~mother
gave her son a sophisticated cell phone capable of sending images and
asked him to continuously transmit her a video of the journey. When the
trip ended, what a misfortune! A left legged creature returned sweet
home!

\noindent
\hskip4pt$\vcenter
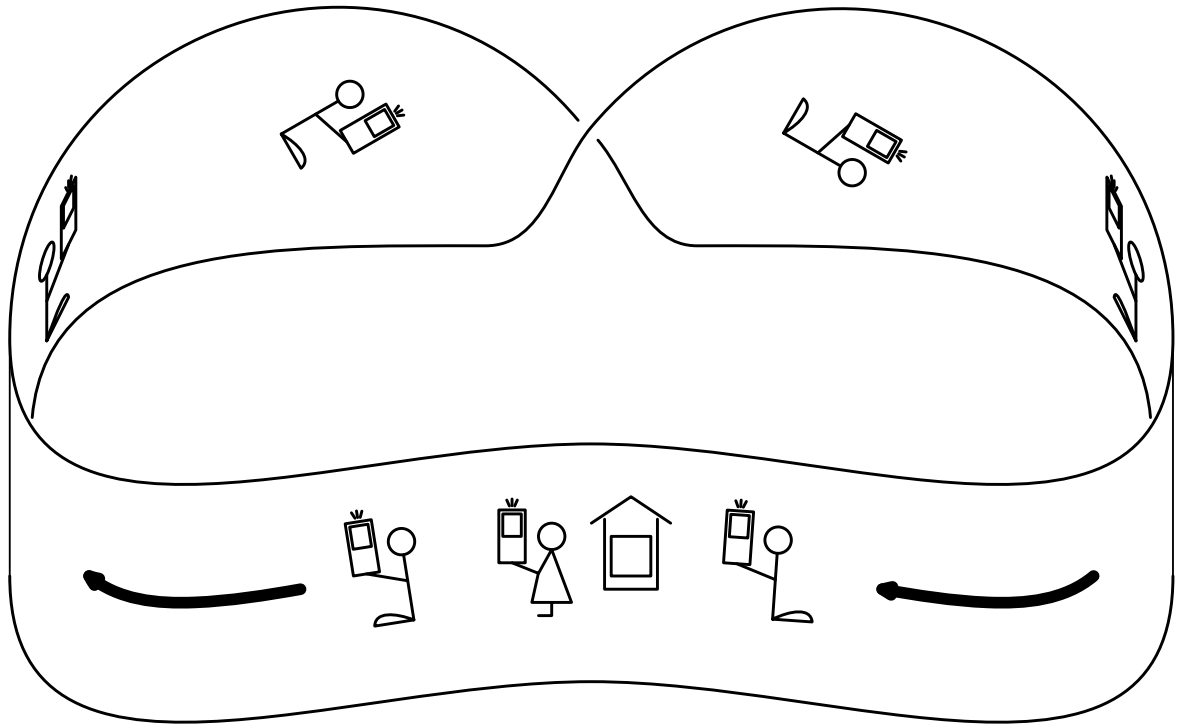
{\hbox{\epsfbox{Pics.4}}}$

\medskip

--- Where is my beloved son? --- asked the desperate mother.

--- And most important $\dots$ where can I buy him a shoe now?

\medskip

Clearly, the last question is merely a Customs and Excise one and can
be solved through an adequate import/export system. More interesting
would be the

\medskip

{\bf1.2.~Question.} At what moment did Candid change his leg?

\medskip

{\bf1.3.~Cogito, ergo sum.} The cartesian coordinates were named after
the French mathematician Ren\'e Descartes. It seems, however, that
Descartes is not to be blamed for disseminating the usage of
coordinates in science. More likely, it was Gottfried Leibniz, one of
Calculus' fathers, the guilty one. Probably, Leibniz also attributed
the name `cartesian coordinates.'

Nobody sees coordinates in Nature. There are no preferred directions
either. (Be careful to apply these ideas in traffic.) In spite of
looking trivial, the above claim has relatively deep consequences.
The~conservation of linear momentum is an example: since there is no
preferred direction, a particle at rest (with respect to some inertial
frame of reference) cannot move spontaneously. In fact, most
conservation laws in physics have a similar origin.\footnote{A rigorous
version of this statement involves the study of symmetries of
differential equations and of the associated conservation laws. Such a
theory was discovered by Emmy Noether, a woman mathematician born in
the city of Erlangen.}

The choice of coordinates while addressing a given problem is
frequently a typical example of an {\it arbitrary choice.} It is not
difficult to realize that an arbitrary choice adds an extra complexity
to the problem. Even worse, such a choice is an obstacle to the
understanding, usually hides subtle features of the problem, and
obscures the essence of the matter.

Every time we are capable of, we are going to avoid arbitrary choices
(of any nature). When an object is essentially related to an arbitrary
choice, we say that it `does not exist.'

\medskip

{\bf1.4.~Guide to the reader.} In what follows, the reader is supposed
either to solve all exercises or to skip (some of) them and accept the
corresponding claims. We have left many hints along the text. There is
also a section entitled `Hints' at the very end. The reader is welcome
to use it from time to time --- once an exercise is solved, one should
look at the corresponding hint anyway. Along the exposition, we use the
exercises as if they were solved.

Some subsections in the book are more advanced and, in principle, may
be not quite `undergraduate.' We believe that the difficulties an
undergraduate student could face in such subsections might be more of
psychological nature than caused by a lack of prerequisites. Anyway,
the more `advanced' subsections are marked with {\bf A}; skipping them
should not compromise the part aimed at undergraduate students.

The book concludes with appendices. They either contain simple and
well-known material (sometimes, in a new exposition) used in the book
or are marked with {\bf A}. Their only common feature is that they are
well inflamed.

\rightheadtext{Projective spaces and their relatives}

\bigskip

\centerline{\bf2.~Projective spaces and their relatives}

\medskip

\vskip10pt

\noindent
\hskip362pt$\vcenter{\hbox{\epsfbox{Pics.3}}}$

\rightskip95pt

\vskip-79pt

In the Euclidean plane $\Bbb E^2$, we fix a point $f$ and consider all
lines passing through~$f$. Such lines constitute points in the space
$\Bbb P_\Bbb R^1$ called the {\it real projective line.} Intuitively,
$\Bbb P_\Bbb R^1$ is one-dimensional. In order to visualize this space,
choose a circle $\Bbb S^1\subset\Bbb E^2$ centred at $f$. The circle
`lists' the lines passing through $f$ : every point $p$ in the circle
generates the line joining $p$ and $f$. Clearly, every line (that is,
every point in $\Bbb P_\Bbb R^1$) is listed exactly twice, by a pair of
diametrically opposed points in the circle. We can therefore visualize
the real projective line as being a `folded'\penalty-10000

\vskip-12pt

\rightskip0pt

\noindent
circle. In this way, we understand that $\Bbb P_\Bbb R^1$ is a circle
itself. The circle $\Bbb P_\Bbb R^1$ can also be obtained from any half
circle contained in $\Bbb S^1$ by simply gluing the ends of the half
circle.

There is another way to visualize $\Bbb P_\Bbb R^1$. We arbitrarily
choose a point in $\Bbb P_\Bbb R^1$ and denote it by $\infty$. This
point corresponds to a line $R_0$ that passes through $f$,
$f\in R_0\subset\Bbb E^2$. We choose a line $T\not\ni f$, parallel to
$R_0$, that does not pass through $f$. The line $T$ will be called the
{\it screen.} Every point $r\in\Bbb P_\Bbb R^1$ (that~is, every line
$R$, $f\in R\subset\Bbb E^2$), except of $\infty$, is displayed on the
screen as the intersection point $R\cap T$. In~this way, the real
projective line is a usual line plus an extra point:
$\Bbb P_\Bbb R^1=\Bbb E^1\sqcup\{\infty\}$, where $\Bbb E^1=T$. We
emphasize again that, {\it a priori,} any point in $\Bbb P_\Bbb R^1$
may play the role of $\infty$.

\medskip

{\bf2.1.~Problem.} Let $R$ be a line in the Euclidean plane $\Bbb E^2$
and let $p$ be a point such that $R\not\ni p\in\Bbb E^2$. Is it
possible, using only a ruler, to construct the line $R'$ passing
through $p$ and parallel to $R$ ?

\medskip

In order to solve Problem 2.1, we need to analyze the concept of
`parallelism' and to discover a `new' mathematical object.

\smallskip

In the Euclidean plane, two distinct lines almost always intersect in a
point. The only exception occurs when the lines are parallel. It would
be nice\footnote{Cicero would say {\it exceptio probat regulam in
casibus non exceptis,} which reads mathematically as `a couple of
counter-examples can substitute a proof of a theorem.' According to
Ivan Karamazov (`The Karamazov brothers' by Fyodor Dostoyevsky)
`$\dots$ they even dare to dream that two parallel lines $\dots$ may
meet somewhere in infinity $\dots$ even if parallel lines do meet and I
see it myself, I shall see it and say that they've met, but still I
won't accept it.'}
if the rule could admit no exception $\dots$

By analogy to the real projective line, we will construct the real
projective plane. In the Euclidean $3$-dimensional space $\Bbb E^3$, we
fix a point $f$ (the light source). The {\it real projective plane\/}
is the set $\Bbb P_\Bbb R^2$ of all lines passing through $f$.

\medskip

{\bf2.2.~Definition.} Let $f\in P\subset\Bbb E^3$ be a plane in
$\Bbb E^3$ passing through $f$. The set $\{R\mid f\in R\subset P\}$
of~all lines $R$ in $P$ passing through $f$ is said to be a {\it
line\/} in $\Bbb P_\Bbb R^2$ ({\it related\/} to $P$). Obviously, this
set is some sort of real projective line $\Bbb P_\Bbb R^1$.

Given two distinct points $r_1,r_2\in\Bbb P_\Bbb R^2$, $r_1\ne r_2$, we
denote by $R_1,R_2\subset\Bbb E^3$ the corresponding lines
in~$\Bbb E^3$. So, there exists a single line in $\Bbb P_\Bbb R^2$ that
`joins' $r_1$ and $r_2$ : the plane $P$ related to the line in question
is the one determined by $R_1,R_2\subset P$. Two distinct lines in
$\Bbb P_\Bbb R^2$ intersect in a single point because the intersection
of two distinct planes that contain $f$ is a line in $\Bbb E^3$ passing
through $f$.

\medskip

\vskip-9pt

{\unitlength=1bp$$\latexpic(350,120)(0,0)\put(0,70){\lline(2,1){75}}
\put(70,0){\lline(2,1){75}}\put(0,70){\lline(1,-1){70}}
\put(75,107.5){\lline(1,-1){70}}\put(69,12){$P$}
\put(34,66){\lline(3,-1){70}}\put(40,66){$R_1$}
\put(43,43){\lline(3,1){65}}\put(90,66){$R_2$}\put(71,58){$f$}
\put(71,45.5){$^{^\bullet}$}\put(200,70){\lline(2,1){62}}
\put(270,0){\lline(2,1){74.5}}\put(200,70){\lline(1,-1){70}}
\put(315,67){\lline(1,-1){29.5}}\put(231,39){\lline(3,1){84}}
\put(271,58){$f$}\put(271,45.5){$^{^\bullet}$}
\put(231,39){\lline(1,2){39}}\put(231,39){\lline(-1,-2){15}}
\put(315,67){\lline(1,2){15}}\put(285,7){\lline(-1,-2){9}}
\put(270.25,117.5){\lline(3,-1){60}}\put(216,9){\lline(3,-1){60}}
\endlatexpic$$}

\vskip-7pt

We will show that the Euclidean plane $\Bbb E^2$ can be seen as a part
of the projective plane $\Bbb P_\Bbb R^2$ in such a way that the lines
in both planes are the `same.'

Indeed, let $f\notin T\subset\Bbb E^3$ be a plane that does not pass
through $f$. Interpreting $f$ as a light source and $T$ as a screen, we
can identify almost every point $r\in\Bbb P_\Bbb R^2$ with its shadow
$p$ on the screen, that is, with the intersection $T\cap R=\{p\}$ of
the screen $T$ with the corresponding line $R\subset\Bbb E^3$. Which
points do not leave a shadow on the screen? Denoting by $P_0$ the plane
that passes through $f$ and is parallel to the screen $T$,
$f\in P_0\subset\Bbb E^3$, we can see that the points that do not have
a shadow on the screen form the line $L_0\simeq\Bbb P_\Bbb R^1$ in
$\Bbb P_\Bbb R^2$ related to $P_0$. In this way, we can see that
$\Bbb P_\Bbb R^2=\Bbb E^2\sqcup\Bbb P_\Bbb R^1$, where $\Bbb E^2=T$ and
$\Bbb P_\Bbb R^1=L_0$.

\vskip3pt

{\unitlength=1bp$$\latexpic(290,103)(0,0)\put(0,0){\lline(0,1){70}}
\put(0,0){\lline(1,1){40}}\put(40,40){\lline(0,1){70}}
\put(0,70){\lline(1,1){40}}\put(26,84){$P_0$}
\put(50,0){\lline(0,1){70}}\put(50,0){\lline(1,1){40}}
\put(90,40){\lline(0,1){70}}\put(50,70){\lline(1,1){40}}
\put(76,84){$T$}\put(4,63){\lline(2,-1){16}}\put(19,47){$^{^\bullet}$}
\put(20,60){$f$}\put(40,45){\lline(2,-1){30}}\put(68,22){$^{^\bullet}$}
\put(69,35){$p$}\put(76.5,26.75){\lline(2,-1){30}}\put(89,8){$R$}
\put(180,0){\lline(0,1){50}}\put(180,60){\lline(0,1){10}}
\put(180,0){\lline(1,1){40}}\put(180,70){\lline(1,1){40}}
\put(220,40){\lline(0,1){70}}\put(206,84){$P_0$}
\put(180,60){\lline(4,-1){40}}\put(184,61){$l'$}
\put(199,47){$^{^\bullet}$}\put(200,60){$f$}
\put(230,0){\lline(0,1){50}}\put(230,60){\lline(0,1){10}}
\put(230,0){\lline(1,1){40}}\put(230,70){\lline(1,1){40}}
\put(270,40){\lline(0,1){70}}\put(256,84){$T$}
\put(230,60){\lline(4,-1){40}}\put(250,57){$l$}
\put(160,60){\lline(1,-1){10}}\put(280,60){\lline(1,-1){10}}
\put(270,50){\lline(-1,0){100}}\put(270,50){\lline(1,0){20}}
\put(180,60){\lline(-1,0){20}}\put(230,60){\lline(-1,0){10}}
\put(270,60){\lline(1,0){10}}\put(224,51.5){$P$}\endlatexpic$$}

\vskip-16pt

{\unitlength=1bp$$\latexpic(170,150)(95,0)\put(42,0){\lline(0,1){26}}
\put(42,60){\lline(0,1){40}}\put(42,0){\lline(1,1){40}}
\put(42,100){\lline(1,1){40}}\put(82,40){\lline(0,1){100}}
\put(68,114){$P_0$}\put(42,60){\lline(4,-1){40}}\put(46,60){$l'$}
\put(61,47){$^{^\bullet}$}\put(62,59){$f$}\put(92,0){\lline(0,1){50}}
\put(92,90){\lline(0,1){10}}\put(92,0){\lline(1,1){40}}
\put(92,100){\lline(1,1){40}}\put(132,40){\lline(0,1){100}}
\put(118,114){$T$}\put(92,90){\lline(4,-1){40}}
\put(132,70){\lline(-4,1){12}}\put(132,50){\lline(-4,1){30}}
\put(22,60){\lline(1,-1){10}}\put(142,60){\lline(1,-1){10}}
\put(132,50){\lline(-1,0){100}}\put(132,50){\lline(1,0){20}}
\put(42,60){\lline(-1,0){20}}\put(132,60){\lline(1,0){10}}
\put(92,90){\lline(-5,-3){10}}\put(132,70){\lline(-5,-2){100}}
\put(132,70){\lline(5,2){30}}\put(162,82){\lline(-1,1){10}}
\put(132,80){\lline(-5,-3){110}}\put(132,80){\lline(5,3){35}}
\put(167,101){\lline(-1,1){21.25}}
\put(145.75,122.25){\lline(-5,-3){13.75}}
\put(32,30){\lline(-1,1){18.6}}\put(13.4,48.6){\lline(5,2){14}}
\put(22,14){\lline(-1,1){21.25}}\put(0.75,35.25){\lline(5,3){16}}
\endlatexpic$$}

\vskip-150pt

\noindent
\hskip265pt$\vcenter
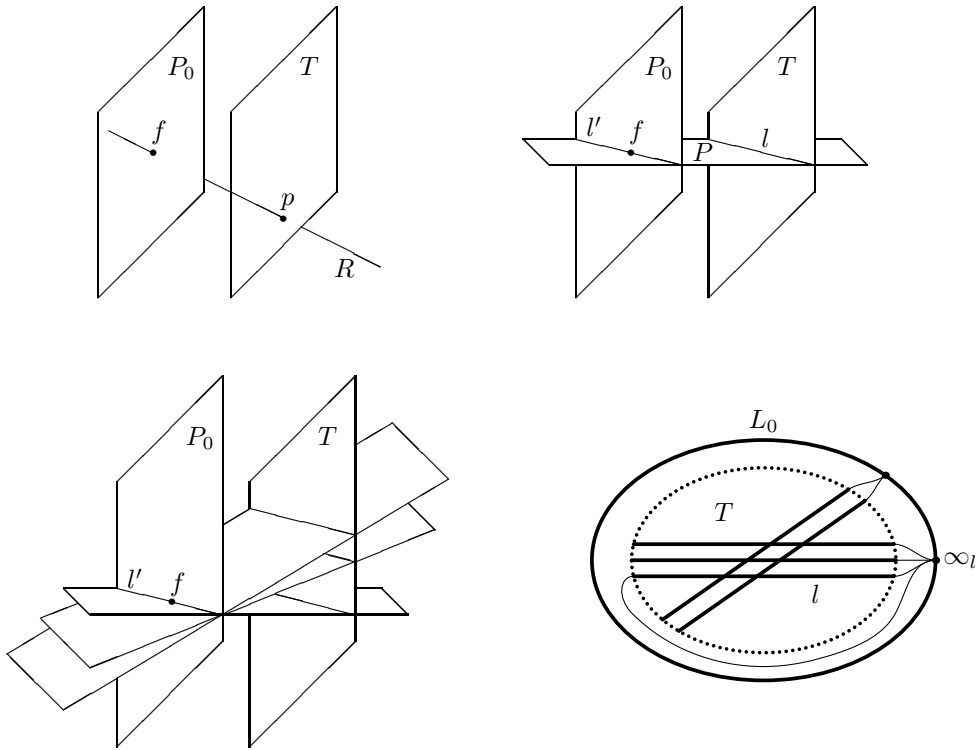
{\hbox{\epsfbox{Pics.1}}}$

\vskip30pt

Let $L\subset\Bbb P_\Bbb R^2$ be a line in $\Bbb P_\Bbb R^2$ distinct
from $L_0$ and let $P$ be the plane related to $L$. So, $P$ is not
parallel to $T$. Therefore, the line $l=T\cap P$ in the plane $T$ is
the shadow of the line $L$ in $\Bbb P_\Bbb R^2$. In the above terms, we
have $L=l\sqcup\{\infty_l\}$, where the point
$\infty_l\in L_0\subset\Bbb P_\Bbb R^2$ corresponds to the line
$l'=P_0\cap P\subset\Bbb E^3$. In this~way, we obtain a one-to-one
correspondence between the lines in $T$ and the lines in
$\Bbb P_\Bbb R^2$ distinct from $L_0$. Each line $l\subset T$ is
extended in $\Bbb P_\Bbb R^2$ by its point at infinity
$\infty_l\in L_0$. The line $L_0\subset\Bbb P_\Bbb R^2$ is formed by
all points at infinity of the lines in $T$.

It is easy to see that two lines are parallel in $T$ iff their points
at infinity are equal. In other words, each family of parallel lines in
$T$ is formed by the lines in $\Bbb P_\Bbb R^2$ that pass through a
same point in $L_0$. Hence, the infinity line $L_0$ can be seen as a
list of such families.

Moving along a line $l\subset T$, independently of the chosen
direction, we finally arrive at the point at infinity
$\infty_l\in L_0$. We arrive at the same point $\infty_l$ if moving
along a line parallel to $l$ and at a different point if moving along a
line non-parallel to $l$.

The following remark is easy, but very important: {\sl Every line in
$\Bbb P_\Bbb R^2$ can be taken as the infinity line.} This solves
Problem 2.1 immediately! Indeed, consider the plane $\Bbb E^2$ as being
inside the real projective plane $\Bbb E^2\subset\Bbb P_\Bbb R^2$ and
use a more powerful ruler that allows us to draw the line in the
projective plane $\Bbb P_\Bbb R^2$ through any two distinct points. Let
us assume that it is possible to construct the parallel line $R'$. Then
we can construct the intersection at infinity $\{q\}=R\cap R'$. Let
$Q$ be the finite set, $p,q\in Q$, of~all points that subsequently
appear during the construction. Such points are intersection points of
lines in $\Bbb P_\Bbb R^2$ that were already constructed at previous
stages plus a finite number of arbitrarily chosen points (that may or
may not belong to the lines that were already constructed). We choose a
new infinity line $L'_0$ in such a way that $L'_0$ passes through no
point in $Q$. We take ${\Bbb E'}^2:=\Bbb P_\Bbb R^2\setminus L'_0$ as a
new (usual) plane. Now, the construction in this new plane
${\Bbb E'}^2$ has to provide the same line $R'$ which, on the other
hand, is not parallel to $R$ because
$\{q\}=R\cap R'\subset{\Bbb E'}^2=\Bbb P_\Bbb R^2\setminus L'_0$. A
contradiction.

\smallskip

Using the more powerful ruler, it is easy to solve the following

\medskip

{\bf2.3.~Exercise.} Let $R_1,R_2$ be distinct parallel lines in the
Euclidean plane $\Bbb E^2$ and let $p\notin R_1,R_2$ be a point,
$R_1,R_2\not\ni p\in\Bbb E^2$. Is it possible, using only a ruler, to
construct the line $R$ passing through $p$ and parallel to $R_1,R_2$ ?

\medskip

Now, we try to visualize the real projective plane $\Bbb P_\Bbb R^2$.
Every sphere $\Bbb S^2\subset\Bbb E^3$ centred at $f$ lists the points
in $\Bbb P_\Bbb R^2$ : each point in $\Bbb P_\Bbb R^2$ is listed twice
by a pair of diametrically opposed points in the sphere. But~this does
not give the faintest idea about the space $\Bbb P_\Bbb R^2$. In order
to understand better the {\it topology\/} of the real projective plane,
we initially cut $\Bbb S^2$ into four pieces and disregard two
redundant ones. Perform-\penalty-10000

\noindent
$\vcenter
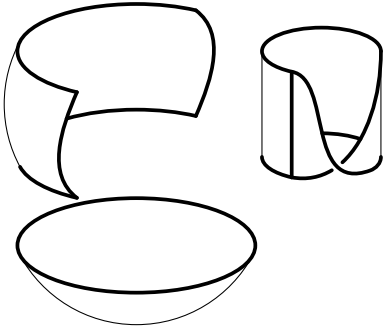
{\hbox{\epsfbox{Pics.2}}}$

\leftskip155pt

\vskip-133pt

\noindent
ing the necessary identifications in one of the two remaining pieces,
we obtain a {\it M\"obius band.} It remains to identify the disc and
the M\"obius band along their boundaries which are circles. In this
way, the structure of the space $\Bbb P_\Bbb R^2$ becomes more or less
clear. Unfortunately, it is impossible to perform such a gluing inside
$\Bbb E^3$.

\medskip

{\bf2.4.~Exercise.} Every line divides the plane $\Bbb E^2$ into two
parts. Into how many parts $4$ generic lines in $\Bbb P_\Bbb R^2$
divide the real projective plane?

\medskip

{\bf2.5.~Exercise.} Visualize the space formed by all unordered pairs
of points in the circle.

\leftskip0pt

\medskip

{\bf2.6.~Projective space.} Let $V$ be a finite-dimensional
$\Bbb K$-linear space, where $\Bbb K=\Bbb R$ or $\Bbb K=\Bbb C$.
We~define the {\it projective space\/} as
$\Bbb P_\Bbb KV:=V^\centerdot/\Bbb K^\centerdot$, where
$V^\centerdot:=V\setminus\{0\}$ is the linear space $V$ `punctured' at
the origin, $\Bbb K^\centerdot$ is (the group of) all non-null elements
in $\Bbb K$, and $V^\centerdot/\Bbb K^\centerdot$ is the {\it
quotient\/} of the {\it action\/} of $\Bbb K^\centerdot$
on~$V^\centerdot$. This means that $V^\centerdot/\Bbb K^\centerdot$ is
the set of equivalence classes in $V^\centerdot$ given by
proportionality with coefficients in $\Bbb K^\centerdot$. (We also
denote $\Bbb P_\Bbb KV=\Bbb P_\Bbb K^n$ if $\dim_\Bbb KV=n+1$.) We have
the quotient map $\pi:V^\centerdot\to\Bbb P_\Bbb KV$ sending every
element to its class. In what follows, we frequently use elements in
$V$ to denote elements in the projective space, that is, we write $p$
in place of $\pi(p)$. In such cases, the~reader is supposed to verify
that our considerations do not change if we rechoose representatives in
$V$ of points in the projective space. One more convention. Given a
subset $S\subset V$, we denote by
$\Bbb P_\Bbb KS:=\pi\big(S\setminus\{0\}\big)\subset\Bbb P_\Bbb KV$ the
image of $S\subset V$ under the quotient map
$\pi:V^\centerdot\to\Bbb P_\Bbb KV$. In this way, for every
$\Bbb K$-linear subspace $K\le V$, we can consider the projective
space $\Bbb P_\Bbb KK$ as a (linear) subspace in~$\Bbb P_\Bbb KV$.

Let $f:V\to\Bbb K$ be a non-null linear functional. We put
$T:=\{v\in V\mid fv=1\}$ and $K:=\ker f$. There is an identification
$(\Bbb P_\Bbb KV\setminus\Bbb P_\Bbb KK)\simeq T$ given by the rule
$v\mapsto\displaystyle\frac v{fv}$. As above, we arrive at the
decomposition
$\Bbb P_\Bbb K^n=\Bbb A_\Bbb K^n\sqcup\Bbb P_\Bbb K^{n-1}$, where the
screen $\Bbb A_\Bbb K^n:=T$ is a $\Bbb K$-affine space (= a
$\Bbb K$-linear space that has forgotten its origin) of dimension $n$.
In terms of $\Bbb P_\Bbb KV$, we can describe $T$ as
$\{p\in\Bbb P_\Bbb KV\mid fp\ne0\}$ and $\Bbb P_\Bbb KK$ as
$\{p\in\Bbb P_\Bbb KV\mid fp=0\}$. (In the expression $fp$, the point
$p$ is to be considered as $p\in V$, but~notice that the equality
$fp=0$ and the inequality $fp\ne0$ do not change their meaning if we
rechoose the representative of the point. --- This is an example of the
above mentioned use of elements in $V$ to denote points in the
projective space.)

Let $x_0,x_1,\dots,x_n:V\to\Bbb K$ be linear coordinates on $V$. We can
define {\it projective coordinates\/} $[x_0,x_1,\dots,x_n]$ on
$\Bbb P_\Bbb KV$ by assuming that
$[kx_0,kx_1,\dots,kx_n]=[x_0,x_1,\dots,x_n]$ for all
$k\in\Bbb K^\centerdot$. Thus, considering each projective coordinate
separately does not provide any meaningful number (but it makes sense
to say whether a coordinate vanishes or not). However, when considered
as entity, the projective coordinates are a mere proportion.

Taking $n+1$ screens, $U_i:=\{p\in\Bbb P_\Bbb KV\mid x_ip\ne0\}$,
$i=0,1,\dots,n$, we have
$\Bbb P_\Bbb KV=\bigcup\limits_{i=0}^nU_i$. Each $U_i$ has $n$ affine
coordinates $y_0,y_1,\dots,y_{i-1},y_{i+1},\dots,y_n$ defined by the
rule $y_j:=x_j/x_i$. The intersection $U_i\cap U_k$ is described as
$U_{ik}:=\big\{p\in U_i\mid y_k(p)\ne0\big\}$ in terms of the
coordinates $y_0,y_1,\dots,y_{i-1},y_{i+1},\dots,y_n$ on $U_i$. The
same intersection is described as
$U_{ki}=\big\{p\in U_k\mid z_i(p)\ne0\big\}$ in terms of the
coordinates $z_0,z_1,\dots,z_{k-1},z_{k+1},\dots,z_n$ on $U_k$. Hence,
$U_{ik}$ is identified with $U_{ki}$ by means of the map
$$U_{ik}\to U_{ki},\qquad(y_0,y_1,\dots,y_{i-1},y_{i+1},\dots,y_n)
\mapsto\Big(\frac{y_0}{y_k},\frac{y_1}{y_k},\dots,\frac{y_{i-1}}{y_k},
\frac{1}{y_k},\frac{y_{i+1}}{y_k},\dots,\frac{y_{k-1}}{y_k},
\frac{y_{k+1}}{y_k},\dots\frac{y_n}{y_k}\Big).$$
In this way, we may interpret $\Bbb P_\Bbb K^n$ as a gluing of $n+1$
copies of $\Bbb A_\Bbb K^n$ identified by the above maps.

For instance, the space $\Bbb P_\Bbb C^1$ can be seen as the gluing of
two copies of $\Bbb C$, equipped with the coordinates~$x_i$, $i=0,1$,
in such a way that the identification between
$U_1\supset U_{10}\simeq\Bbb C^\centerdot$ and
$U_0\supset U_{01}\simeq\Bbb C^\centerdot$ is given by the formula
$x_0x_1=1$. In particular, we visualize $\Bbb P_\Bbb C^1$ as
$\big\{[1,x_1]\mid x_1\in\Bbb C\big\}\simeq\Bbb C$ extended by the
point at infinity $\infty=[0,1]$.

\medskip

\vskip5pt

\noindent
\hskip312pt$\vcenter{\hbox{\epsfbox{Pics.5}}}$

\rightskip145pt

\vskip-132pt

{\bf2.7.~Sphere and stereographic projection.} Let $V$ be an
$\Bbb R$-linear space, $\dim_\Bbb RV=n+1$. We define the
{\it$n$-sphere\/} as being $\Bbb S^n:=V^\centerdot/\Bbb R^+$, where
$R^+:=\{r\in\Bbb R\mid r>0\}$. A more common definition of the unit
$n$-sphere inside Euclidean space is
$\Bbb S^n:=\big\{p\in\Bbb E^{n+1}\mid\allowmathbreak\langle
p,p\rangle=\nomathbreak1\big\}$,
where $\langle-,-\rangle$ stands for the usual inner product in
$\Bbb E^{n+1}$. We define the {\it tangent space\/} $\T_p\Bbb S^n$ to
$\Bbb S^n$ at $p\in\Bbb S^n$ as
$\T_p\Bbb S^n:=\allowmathbreak p^\perp\le\nomathbreak\Bbb E^{n+1}$. In
order to have a hyperplane that is indeed tangent to the sphere at $p$,
it is better to take $p+p^\perp$ in place of $p^\perp$, but we prefer
the above definition as it provides an obvious linear space. The {\it
stereographic projection\/}
$\varsigma_p:\Bbb S^n\setminus\{-p\}\to\T_p\Bbb S^n$ sends the point
$q\in\Bbb S^n\setminus\{-p\}$ to the intersection
$\T_p\Bbb S^n\cap R(-p,q)$, where $R(-p,q)$ denotes the line joining
$-p$ and $q$.

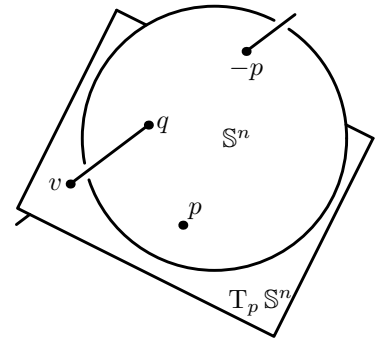
\rightskip0pt

When $n=2$, we can interpret the stereographic projection as
`unwrapping' the sphere punctured at the point $-p$ into the plane
tangent to the sphere at the point $p$. This unwrapping is one of the
typical ways of exhibiting geographic maps.

\medskip

{\bf2.8.~Exercise.} Prove the explicit formulae

\smallskip

\centerline{$\varsigma_p:\Bbb S^n\setminus\{-p\}\ni
q\mapsto\displaystyle\frac{q+p}{1+\langle
q,p\rangle}-p\in{\T}_p\Bbb S^n,\qquad\varsigma_p^{-1}:{\T}_p\Bbb S^n\ni
v\mapsto\frac{2(v+p)}{1+\langle v,v\rangle}-p\in\Bbb
S^n\setminus\{-p\}$.}

\medskip

\vskip5pt

\noindent
\hskip295pt$\vcenter{\hbox{\epsfbox{Pics.6}}}$

\rightskip165pt

\vskip-140pt

{\bf2.9.~Riemann sphere.} Using a couple of stereographic projections,
we can see that $\Bbb P_\Bbb C^1\simeq\Bbb S^2$. Indeed, let us treat
the tangent planes $\T_p\Bbb S^2$ and $\T_{-p}\Bbb S^2$ to the unit
sphere $\Bbb S^2$ at the diametrically opposed points
$p,-p\in\Bbb S^2\subset\Bbb E^3$ as being planes of complex numbers,
$\T_p\Bbb S^2\simeq\Bbb C_0$ and $\T_{-p}\Bbb S^2\simeq\Bbb C_1$, in
such a way that the real axes are parallel with the same directions and
the imaginary axes are parallel with the opposite directions. (In order
to facilitate the visualization, we draw the tangent planes as passing
through the points, $p\in\T_p\Bbb S^2$ and $-p\in\T_{-p}\Bbb S^2$.)
Let~$0\ne x\in\Bbb C_0\simeq\T_p\Bbb S^2=p^\perp$. Applying the
formulae from Exercise 2.8, we obtain
$\varsigma_{-p}\varsigma_p^{-1}x=x/\langle x,x\rangle=1/\overline x$,
which corresponds to $1/x\in\Bbb C_1\simeq\T_{-p}\Bbb S^2$. In other
words, the gluing of the planes\break

\vskip-12pt

\rightskip0pt

\noindent
$\T_p\Bbb S^2$ and $\T_{-p}\Bbb S^2$ resulting in $\Bbb S^2$ is the
same as the above described gluing of $U_0$ and $U_1$ resulting
in $\Bbb P_\Bbb C^1$.

\medskip

{\bf2.10.~Exercise.} Prove that the stereographic projection
$\varsigma_p$ establishes a one-to-one correspondence between
subspheres in $\Bbb S^n$ (= intersections of $\Bbb S^n$ with affine
subspaces in $\Bbb E^{n+1}$) and subspheres or affine subspaces in
$\T_p\Bbb S^n$.

\medskip

{\bf2.11.~Exercise.} Prove that the stereographic projection preserves
angles between curves.

\medskip

{\bf2.12A.~Grassmannians.} We take and fix finite-dimensional
$\Bbb K$-linear spaces $P,V$ and denote by
$$M:=\big\{p\in{\Lin}_\Bbb K(P,V)\mid\ker p=0\big\}$$
the open subset of all monomorphisms in the $\Bbb K$-linear space
$\Lin_\Bbb K(P,V)$. The group $\GL_\Bbb KP$ of all nondegenerate
$\Bbb K$-linear transformations of $P$ acts from the right on
$\Lin_\Bbb K(P,V)$ and on $M$. By definition, the {\it grassmannian\/}
$\Gr_\Bbb K(k,V)$ is the quotient space
$${\Gr}_\Bbb K(k,V):=M/{\GL}_\Bbb KP,\qquad\pi:M\to M/{\GL}_\Bbb KP,$$
where $k:=\dim_\Bbb KP$. It is the space of all $k$-dimensional
$\Bbb K$-linear subspaces in $V$. In the case of $\Bbb K=\Bbb R$, we
can also take the group
$\GL_\Bbb R^+P:=\{g\in\GL_\Bbb RP\mid\det g>0\}$ in place of
$\GL_\Bbb KP$, obtaining the {\it grassmannian\/}
$${\Gr}_\Bbb R^+(k,V):=M/{\GL}_\Bbb R^+P,\qquad\pi':M\to M/{\GL}_\Bbb
R^+P$$
of {\it oriented\/} $k$-dimensional $\Bbb R$-linear subspaces in $V$.

\rightheadtext{Smooth spaces and smooth functions}

\bigskip

\centerline{\bf3.~Smooth spaces and smooth functions}

\medskip

Why do we feel that the $2$-sphere is smooth and the (surface of a
$3$-) cube is not? We guess that the concept of smooth function answers
well this question. Everybody knows, at least at the level of
intuition, what a smooth function is.\footnote{The following story
about the `Grothendieck prime' (Alexander Grothendieck, one of the
greatest mathematicians of our times) comes to mind. Somebody
suggested: `Pick a prime number.' Grothendieck replied: `You mean like
$57$ ?'}
Actually, instead of any kind of formal definition, it seems better to
simply list the properties of (smooth) functions that we are going to
use. Inevitably, we are to simultaneously introduce the properties of
(smooth) spaces.

In this section, we try to focus ourselves on understanding and
clarifying the nature of objects and concepts. It turns out that our
introduction to differential topology came out a little bit
nonstandard, but this pays off: the same exposition works for
algebraic/complex geometry. The reader is welcome to get back to this
material and give it a broader look; however, in the first reading, one
may opt to stuck with usual smooth functions and spaces.

\medskip

{\bf3.1.~Introductory remarks.} We fix some field $\Bbb K$. In our
applications, it will be the field $\Bbb R$ of real numbers or the
field $\Bbb C$ of complex numbers. We would like to speak of local
$\Bbb K$-valued `smooth' functions
$M\opensup U\overset f\to\longrightarrow\Bbb K$ defined on open subsets
$U\opensub M$ of a given topological space $M$. Denote by $\Cal F$ all
such functions and by $\Cal F(U)$, those with a given $U\opensub M$. We
can sum and multiply the functions in $\Cal F(U)$. Naturally, the
constant functions should be included in $\Cal F(U)$. In other words,
$\Cal F(U)$ is a commutative $\Bbb K$-algebra. A more important feature
of `smooth' functions is that this concept is {\bf local}. This means
that, for  $W\opensub U\opensub M$ and $f\in\Cal F(U)$, the restriction
$f|_W:W\to\Bbb K$ belongs to $\Cal F(W)$ and {\it vice versa\/}: if a
function is locally `smooth,' it must be `smooth.' Thus, we arrive at
the following definition.

\medskip

{\bf3.2.~Sheaves of functions.} Let $M$ be a topological space and let
$\Cal F:=\bigsqcup\limits_{U\opensub M}\Cal F(U)$ be a collection of
$\Bbb K$-valued functions such that $\Cal F(U)$ is a $\Bbb K$-algebra
for every $U\opensub M$ and the following conditions hold.

\smallskip

\leftskip60pt\rightskip60pt

\noindent
$\bullet$ If $W\opensub U\opensub M$ and $f\in\Cal F(U)$, then
$f|_W\in\Cal F(W)$.

\noindent
$\bullet$ Let us be given open subsets $U_i\opensub M$, $i\in I$, and a
function $U\overset f\to\longrightarrow\Bbb K$, where
$U:=\bigcup\limits_{i\in I}U_i$. If $f|_{U_i}\in\Cal F(U_i)$ for every
$i\in I$, then $f\in\Cal F(U)$.

\leftskip0pt\rightskip0pt

\smallskip

\noindent
Then $\Cal F$ is a {\it sheaf\/} of $\Bbb K$-valued functions on $M$.

Speaking slightly informally, a sheaf of functions corresponds to a
local property of a $\Bbb K$-valued function preserved by the
$\Bbb K$-algebra operations.

{\unitlength=1bp$$\latexpic(57,57)(-170,0)
\put(0,60){$M\opensup U\opensup W\ni p$}\put(4,42){$\Cal F(U)$}
\put(29,45){\vector(1,0){20}}\put(34,35){$|_W$}\put(52,42){$\Cal F(W)$}
\put(20,37){\vector(1,-2){14}}\put(54,37){\vector(-1,-2){14}}
\put(33,0){$\Cal F_p$}
\endlatexpic$$}

\rightskip94pt

\vskip-87pt

Let $p\in M$ be fixed, let $p\in U_1,U_2\opensub M$, and let
$f_i\in\Cal F(U_i)$, $i=1,2$. We write $f_1\sim f_2$ if there exists
$U\opensub U_1\cap U_2$ such that $p\in U$ and $f_1|_U=f_2|_U$.
Obviously, $\sim$ is an equivalence relation. The~corresponding
equivalence class $f_p$ is the {\it germ\/} of $f\in\Cal F$ at $p$. All
germs at $p$ form the {\it stalk\/} $\Cal F_p$ of $\Cal F$ at $p$. The
stalk is a $\Bbb K$-algebra and, for $p\in U\opensub M$, we have the
homomorphism $\Cal F(U)\to\Cal F_p$, $f\mapsto f_p$,
of~$\Bbb K$-algebras which is compatible with restrictions.

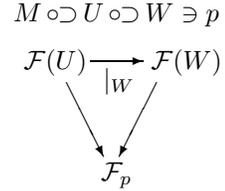
\rightskip0pt

The $\Bbb K$-algebra $\Cal F_p$ splits into $\Bbb K$ (the constants)
and the ideal
$\frak m_p:=\big\{f_p\mid f(p)=0\big\}\triangleleft\Cal F_p$ formed by
the germs that vanish at $p$. So, $\Cal F_p=\Bbb K+\frak m_p$.

\medskip

{\bf3.3.~Basic example.} Let $V$ be a finite-dimensional
$\Bbb K$-linear space equipped with the usual topology. Let
$p\in U\opensub V$, $f:U\to\Bbb K$, and $v\in V$ be a point, a
function, and a vector. We denote by
$$v_pf:=
\lim_{\varepsilon\to0}\frac{f(p+\varepsilon v)-f(p)}\varepsilon$$
the {\it$v$-directional derivative\/} of $f$ at $p$. If
$f:=\varphi|_U$, where $\varphi\in V^*:=\Lin_\Bbb K(V,\Bbb K)$ is a
$\Bbb K$-linear functional, then such a derivative exists and equals
$v_pf=\varphi v$. Of course, $v_pc=0$ for any constant function $c$.
If~$v_pf$ exists for every $p\in U$, we define the {\it partial
derivative\/} $[v]_Uf:U\mapsto\Bbb K$ by the rule
$[v]_Uf:p\mapsto v_pf$. A~continuous function $f:U\to\Bbb K$ is said to
be {\it smooth\/} of {\it class\/} $C^0$. By induction, a function
$f:U\to\Bbb K$ is {\it smooth\/} of {\it class\/} $C^k$ iff the
function $[v]_Uf:U\to\Bbb K$ (exists and) is smooth of class $C^{k-1}$
for every $v\in V$. A function $f:U\to\Bbb K$ is {\it smooth\/} (of
{\it class\/} $C^\infty$) iff it is smooth of class $C^k$ for every
$k\ge0$.

\smallskip

{\bf3.3.1.~Exercise.} Let $f_1,f_2:U\to\Bbb K$, $p\in U\opensub V$, and
$v\in V$ be such that $v_pf_1,v_pf_2$ exist. Show that
$v_p(f_1+f_2),v_p(f_1f_2)$ exist and
$$v_p(f_1+f_2)=v_pf_1+v_pf_2,\qquad
v_p(f_1f_2)=f_1(p)v_pf_2+f_2(p)v_pf_1.$$
(The latter is the well-known {\it Leibniz rule.}) Show that $C^k$,
formed by all smooth functions of class $C^k$, $0\le k\le\infty$, is a
sheaf of $\Bbb K$-valued functions on $V$. We have
$C^k(U)\subset C^{k-1}(U)$ and $[v]_U:C^k(U)\to C^{k-1}(U)$ for all
$v\in V$ and $U\opensub V$. Note that $[v]_U$ is compatible with
restrictions. So, we can write $[v]$ instead of $[v]_U$.

\smallskip

{\bf3.3.2.~Exercise.} If $v_pf$ exists, then $(kv)_pf$ exists and
$(kv)_pf=kv_pf$ for every $k\in\Bbb K$. For $f\in C^1$ and
$v,w\in V$, we have $[v+w]f=[v]f+[w]f$.

\smallskip

{\bf3.3.3.~Exercise {\rm(Taylor's formula)}.} Let $p\in V$ and
$g\in C_p^\infty$. Then there exist a unique linear functional
$\varphi\in V^*$ and $h\in\frak m_p^2$ such that
$g=g(p)+\varphi_p-\varphi p+h$.

\smallskip

{\bf3.3.4.~Exercise.} Show that the topology on $V$ is the weakest one
such that all functions $C^\infty(V)\ni f:V\to\Bbb K$ are continuous.

\smallskip

{\bf3.3.5.~Exercise.} Let $V\opensup U\overset\psi\to\longrightarrow W$
be a map into a finite-dimensional $\Bbb K$-linear space $W$. Suppose
that $W^*\circ\psi\subset C^\infty(U)$. Show that $\psi$ is continuous
and that $f\circ\psi\in C^\infty\big(\psi^{-1}(X)\big)$ for all
$X\opensub W$ and $f\in C^\infty(X)$.

\smallskip

Until the end of this section, the reader may assume for simplicity
that the sheaves we deal with are all induced by the sheaves
$C^\infty$.

\medskip

{\bf3.4.~Smooth maps and induced structures.} Let $(M_1,\Cal F_1)$ and
$(M_2,\Cal F_2)$ be spaces with sheaves of functions. A continuous map
$\psi:M_1\to M_2$ is `{\it smooth\/}' if
$f_2\circ\psi\in\Cal F_1\big(\psi^{-1}(U_2)\big)$ for all
$U_2\opensub M_2$ and $f_2\in\Cal F_2(U_2)$.

Let $(M_2,\Cal F_2)$ be a space with a sheaf of functions and let
$\varphi:M\to M_2$ be a map. Then there exist a weakest topology and a
smallest sheaf $\Cal F$ of functions on $M$ such that $\varphi$ is
smooth. More precisely, the~open subsets in $M$ are of the form
$U=\varphi^{-1}(U_2)$, where $U_2\opensub M_2$. A function
$M\opensup U\overset f\to\longrightarrow\Bbb K$ belongs to $\Cal F(U)$
iff it is locally of the form $f_2\circ\varphi$, i.e., iff there exist
an open cover $U_2=\bigcup\limits_{i\in I}U_i$ and functions
$f_i\in\Cal F_2(U_i)$ such that $U=\varphi^{-1}(U_2)$ and
$f|_{\varphi^{-1}(U_i)}=f_i\circ\varphi$ for all $i\in I$. The
introduced structure on $M$\break

{\unitlength=1bp$$\latexpic(30,33)(-183,0)
\put(0,41){$M_1$}\put(14,44){\vector(1,0){22}}\put(22,35){$\psi$}
\put(39,41){$M_2$}\put(10,37){\vector(1,-2){14}}
\put(11,15){$\vartheta$}\put(28,9){\vector(1,2){14}}
\put(37,19){$\varphi$}\put(21,0){$M$}
\endlatexpic$$}

\rightskip63pt

\vskip-74pt

\noindent
is called {\it induced\/} by $\varphi$. It is universal in the
following sense. If $\psi=\varphi\circ\vartheta$ for some map
$\vartheta:M_1\to M$ and a smooth map
$(M_1,\Cal F_1)\overset\psi\to\longrightarrow(M_2,\Cal F_2)$, then
$\vartheta$ is smooth. The~concept of induced structure usually applies
to subsets $M\subset M_2$. In this case, the induced sheaf is denoted
by $\Cal F_2|_M$. In the easy (and important) case of $M\opensub M_2$,
we have $\Cal F_2|_M=\bigsqcup\limits_{U\opensub M}\Cal F_2(U)$.

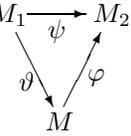
\rightskip0pt

Let $(M_1,\Cal F_1)$ be a space with a sheaf of functions and let
$\varphi:M_1\to M$ be a map. Then there exist a strongest topology and
a largest sheaf $\Cal F$ of functions on $M$ such that $\varphi$ is
smooth. More precisely, $U\opensub M$ iff $\varphi^{-1}(U)\opensub M_1$
and $M\opensup U\overset f\to\longrightarrow\Bbb K$ belongs to
$\Cal F(U)$ iff $f\circ\varphi\in\Cal F_1\big(\varphi^{-1}(U)\big)$.
The introduced\break

{\unitlength=1bp$$\latexpic(30,30)(-183,0)
\put(0,41){$M_1$}\put(14,44){\vector(1,0){22}}\put(22,35){$\psi$}
\put(39,41){$M_2$}\put(10,37){\vector(1,-2){14}}\put(10,18){$\varphi$}
\put(28,9){\vector(1,2){14}}\put(37,19){$\vartheta$}\put(21,0){$M$}
\endlatexpic$$}

\rightskip66pt

\vskip-72pt

\noindent
structure on $M$ is called the {\it quotient\/} by $\varphi$. It is
universal in the following sense. If~$\psi=\vartheta\circ\varphi$ for
some map $\vartheta:M\to M_2$ and a smooth map
$(M_1,\Cal F_1)\overset\psi\to\longrightarrow(M_2,\Cal F_2)$, then
$\vartheta$ is smooth. The concept of quotient structure usually
applies to the quotient by an equivalence relation
$M_1\to M:=M_1/\sim$.

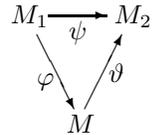
\rightskip0pt

\smallskip

{\bf3.4.1.~Exercise.} Let $(M,\Cal F)$ and $(N,\Cal G)$ be spaces with
sheaves of functions, let $\psi:M\to N$ be a~map, and let
$N=\bigcup\limits_{i\in I}U_i$ and
$\psi^{-1}(U_i)=\bigcup\limits_{j\in J_i}U_{ij}$, $i\in I$, be open
covers. Show that $\psi$ is smooth iff all
$\psi|_{U_{ij}}:U_{ij}\to U_i$ are smooth, where $U_i$ and $U_{ij}$ are
equipped with the induced structures. In other words, the concept of a
smooth map is local.

\smallskip

{\bf3.4.2.~Exercise.} Let $M$ be a set and suppose that
$M=\bigcup\limits_{i\in I}M_i$, where every $M_i$ is equipped with a
topology and a sheaf $\Cal F_i$ of $\Bbb K$-valued functions such that
$M_i\cap M_j\opensub M_i$ and
$\Cal F_i|_{M_i\cap M_j}=\Cal F_j|_{M_i\cap M_j}$ for all $i,j\in I$.
Verify that there exist a unique topology and a sheaf $\Cal F$ on $M$
such that $M_i\opensub M$ and the structure on $M_i$ is induced by that
on $M$ for all $i\in I$. In this situation, we say that $(M,\Cal F)$ is
a {\it gluing\/} of $(M_i,\Cal F_i)$, $i\in I$. We have already seen a
couple of examples of gluing in Subsections 2.6 and 2.9.

\smallskip

{\bf3.4.3.~Example.} Let $V$ be an $\Bbb R$-linear space with
$\dim_\Bbb RV=n+1$. Then $V^\centerdot:=V\setminus\{0\}\opensub V$ gets
the induced $C^\infty$-structure. If $v_2=rv_1$ for some $r>0$, we
write $v_1\sim v_2$. Then we obtain the quotient structure on the
$n$-sphere $\Bbb S^n:=V^\centerdot/\sim$ and the smooth map
$\pi:V^\centerdot\to\Bbb S^n$.

\smallskip

{\bf3.4.4A.~Example.} More generally, the grassmannians
$\pi:M\to\Gr_\Bbb K(k,V)$ and $\pi':M\to\Gr_\Bbb R^+(k,V)$ (see 2.12A)
are equipped with the quotient structure.

\smallskip

{\bf3.4.5.~Example.} Let $V$ be an Euclidean $\Bbb R$-linear space with
$\dim_\Bbb RV=n+1$. Then
$S:=\big\{v\in V\mid\langle v,v\rangle=1\big\}\subset V$ is closed. We
have the induced $C^\infty$-structure on $S\subset V^\centerdot$.

\smallskip

{\bf3.4.6.~Exercise.} Show that the composition
$S\hookrightarrow V^\centerdot\overset\pi\to\longrightarrow\Bbb S^n$
(see Examples 3.4.3 and 3.4.5) is a diffeomorphism (i.e., a smooth
isomorphism).

\medskip

{\bf3.5.~Product and fibre product.} We fix a certain class $\Cal C$ of
spaces with sheaves of $\Bbb K$-valued functions and assume that
$\Cal C$ is closed with respect to taking open subspaces (equipped with
the induced structure) and with respect to gluing. So, for a gluing
$M=\bigcup\limits_{i\in I}M_i$, we have $M_i\in\Cal C$ for all $i\in I$
iff $M\in\Cal C$ (actually, we will need only the gluings with
countable or finite $I$). In other words, the property `to belong to
$\Cal C$' is local.

\smallskip

{\unitlength=1bp$$\latexpic(40,32)(-130,0)
\put(50,33){$M$}\put(47,32){\vector(-3,-2){34}}\put(23,26){$\psi_1$}
\put(62,32){\vector(3,-2){34}}\put(80,24){$\psi_2$}
\put(55,30){\vector(0,-1){20}}\put(57,19){$\psi$}\put(0,0){$M_1$}
\put(34,2){\vector(-1,0){19}}\put(23,5){$\pi_1$}
\put(36,0){$M_1\times M_2$}\put(77,2){\vector(1,0){19}}
\put(80,5){$\pi_2$}\put(98,0){$M_2$}
\endlatexpic$$}

\rightskip125pt

\vskip-63pt

{\bf3.5.1.~Product.} Let $M_1,M_2\in\Cal C$. A structure on
$M_1\times M_2$ providing $M_1\times M_2\in\Cal C$ is a $\Cal C$-{\it
product\/} if the projections $\pi_i:M_1\times M_2\to M_i$ are smooth
and, for any $M\in\Cal C$ and smooth maps $\psi_i:M\to M_i$, the map
$\psi:M\to M_1\times M_2$ in the commutative diagram is smooth.

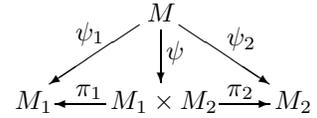
\rightskip0pt

\smallskip

{\unitlength=1bp$$\latexpic(40,65)(-123,0)
\put(18,66){$M_1\times M_2$}\put(42,63){\vector(1,-1){21}}
\put(38,63){\vector(0,-1){54}}\put(32,43){$1$}
\put(34,63){\vector(-1,-1){21}}\put(18,0){$M_1\times M_2$}
\put(49,33){$M_1\times'M_2$}\put(47,36){\lline(-1,0){6}}
\put(64,30){\vector(-1,-1){21}}\put(32,10){\vector(-1,1){19}}
\put(60,63){\vector(2,-1){42}}\put(60,9){\vector(2,1){42}}
\put(0,33){$M_1$}\put(35,36){\vector(-1,0){19}}
\put(93,36){\vector(1,0){10}}\put(104,33){$M_2$}
\endlatexpic$$}

\rightskip133pt

\vskip-97pt

{\bf3.5.2.~Exercise.} Let $M_1,M_2\in\Cal C$. Show that a
$\Cal C$-product structure on $M_1\times M_2$ is unique if it exists.

\smallskip

{\bf3.5.3.~Exercise.} Let
$S_1,S_2,S_1\times S_2,M_1,M_2,M_1\times M_2\in\Cal C$, where\break

{\unitlength=1bp$$\latexpic(40,47)(203,0)
\put(45,51){$M$}\put(43,51){\vector(-2,-1){30}}
\put(51,48){\vector(0,-1){11}}\put(56,51){\vector(2,-1){30}}
\put(2,27){$S_1$}\put(32,30){\vector(-1,0){18}}
\put(34,27){$S_1\times S_2$}\put(68,30){\vector(1,0){19}}
\put(89,27){$S_2$}\put(7,24){\vector(0,-1){15}}
\put(51,24){\vector(0,-1){15}}\put(94,24){\vector(0,-1){15}}
\put(0,0){$M_1$}\put(30,3){\vector(-1,0){15}}
\put(31,0){$M_1\times M_2$}\put(72,3){\vector(1,0){15}}
\put(88,0){$M_2$}
\endlatexpic$$}

\leftskip113pt

\vskip-91pt

\noindent
$S_i\subset M_i$, $i=1,2$, and $S_1\times S_2\subset M_1\times M_2$ are equipped with the
induced structures and $M_1\times M_2$ is a $\Cal C$-product. Prove
that $S_1\times S_2$ is a $\Cal C$-product.

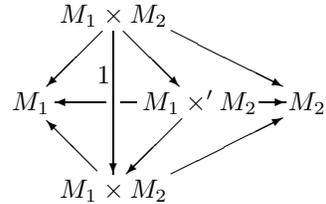
\rightskip0pt

\smallskip

{\bf3.5.4.~Exercise.} Let $M_i,N_j\in\Cal C$ for all $i\in I$ and
$j\in J$ and let $M=\bigcup\limits_{i\in I}M_i$\break

\vskip-5pt

\rightskip163pt

{\unitlength=1bp$$\latexpic(40,58)(-90,0)
\put(13,58){$M_i$}\put(88,61){\vector(-1,0){60}}
\put(90,58){$M_i\times N$}\put(123,57){\vector(1,-1){19}}
\put(18,39){\vector(0,1){15}}\put(109,39){\vector(0,1){15}}
\put(0,29){$M_i\cap M_j$}\put(56,32){\vector(-1,0){15}}
\put(58,29){$(M_i\cap M_j)\times N$}\put(125,32){\vector(1,0){15}}
\put(143,29){$N$}\put(13,0){$M_j$}\put(88,3){\vector(-1,0){60}}
\put(90,0){$M_j\times N$}\put(125,9){\vector(1,1){17}}
\put(18,26){\vector(0,-1){15}}\put(109,26){\vector(0,-1){15}}
\endlatexpic$$}

\vskip-92pt

\noindent
and $N=\bigcup\limits_{j\in J}N_j$ be gluings. Suppose\break

\vskip-5pt

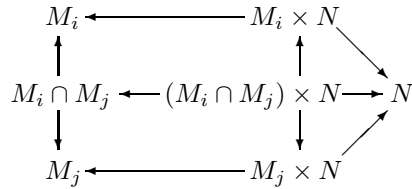
\leftskip0pt

\noindent
that there exists a $\Cal C$-product structure on $M_i\times N_j$ for
all $i\in I$ and $j\in J$. Show that the gluing of $M_i\times N_j$
provides a $\Cal C$-product structure on $M\times N$.

\smallskip

{\bf3.5.5.~Fibre product.} Let $M_1,M_2,B\in\Cal C$ and let
$\varphi_i:M_i\to B$ be smooth maps, $i=1,2$. We define\newline

\vskip-5pt

\rightskip0pt

$$M_1\times_BM_2:=\big\{(p_1,p_2)\in M_1\times
M_2\mid\varphi_1(p_1)=\varphi_2(p_2)\big\},$$
$$\pi_i:M_1\times_BM_2\to M_i,\quad\pi_i:(p_1,p_2)\mapsto p_i,\quad
i=1,2.$$

{\unitlength=1bp$$\latexpic(40,75)(-123,0)
\put(53,74){$M$}\put(51,73){\vector(-3,-2){40}}\put(21,63){$\psi_1$}
\put(64,73){\vector(3,-2){40}}\put(88,61){$\psi_2$}
\put(58,71){\vector(0,-1){24}}\put(60,56){$\psi$}\put(0,37){$M_1$}
\put(34,39){\vector(-1,0){19}}\put(23,42){$\pi_1$}
\put(36,37){$M_1\times_BM_2$}\put(84,39){\vector(1,0){19}}
\put(87,42){$\pi_2$}\put(105,37){$M_2$}\put(11,33){\vector(3,-2){40}}
\put(35,21){$\varphi_1$}\put(104,34){\vector(-3,-2){41}}
\put(74,23){$\varphi_2$}\put(54,0){$B$}
\endlatexpic$$}

\rightskip131pt

\vskip-107pt

\noindent
Clearly, $\varphi_1\circ\pi_1=\varphi_2\circ\pi_2$. A structure on
$M_1\times_BM_2$ providing $M_1\times_BM_2\allowmathbreak\in\Cal C$ is
a {\it fibre product\/} in $\Cal C$ (or a {\it fibre product\/}
$\Cal C$-structure) if $\pi_1,\pi_2$ are smooth and, for any
$M\in\Cal C$ and smooth maps
$M_1\overset\psi_1\to\longleftarrow
M\overset\psi_2\to\longrightarrow M_2$
such that $\varphi_1\circ\psi_1=\varphi_2\circ\psi_2$, the map
$\psi:M\to M_1\times_BM_2$ in the commutative diagram is smooth.

It is frequently useful to visualize the fibre product $M_1\times_BM_2$
as a family of products parameterized by $B$. More specifically,
$M_1\times_BM_2=$\break

\vskip-12pt

\rightskip0pt

\noindent
$\bigsqcup\limits_{p\in B}\varphi_1^{-1}(p)\times\varphi_2^{-1}(p)$,
where $\varphi_1^{-1}(p)\times\varphi_2^{-1}(p)$ is the product of the
{\it fibres\/} of $\varphi_1$ and $\varphi_2$ over $p\in B$.

\smallskip

{\bf3.5.6.~Exercise.} Let $M_1,M_2,B\in\Cal C$. Show that a fibre
product $\Cal C$-structure on $M_1\times_BM_2$ is unique if it exists.

\smallskip

{\bf3.5.7.~Exercise.} Let
$M_1,M_2,B,M_1\times M_2,M_1\times_BM_2\in\Cal C$, where
$M_1\times M_2$ is a $\Cal C$-product and
$M_1\times_BM_2\subset M_1\times M_2$ is equipped with the induced
structure. Prove that $M_1\times_BM_2$ is a fibre product in $\Cal C$.

\medskip

{\bf3.6.~Tangent bundle.} We need to understand what is a tangent
vector at a point $p\in M$ to a space $M$ equipped with a sheaf of
functions. Everybody seems to `know' what a tangent vector to a smooth
surface $M\subset\Bbb K^3$ is and can even draw it when
$\Bbb K=\Bbb R$. Nevertheless, there are a couple of problems. The
first consists in the words `smooth surface' --- we did not yet define
a smooth subspace and the definition that first comes to mind tends to
use the concept of a tangent vector itself $\dots$ The other problem is
even more heavy. Our intuitive view on a tangent vector is in no way
intrinsic. So, we have no clear idea on how to compare tangent vectors
at the same point $p\in M$ that come from different smooth embeddings
$M\hookrightarrow\Bbb K^n$.

Fortunately, both problems can be solved with the same remedy. For the
first, we can restrict the sheaf $\Cal F$ on $\Bbb K^n$ to $M$ and hope
to characterize the smoothness of $M$ in terms of $\Cal F|_M$. Our
basic example~3.3 provides a hint on how to manage the second problem.
We can simply interpret an intuitive tangent vector $v$ at $p\in M$ as
being a derivative in its direction. It is true that the expression
$f(p+\varepsilon v)$ makes no\break

\vskip-5pt

\noindent
\hskip324pt$\vcenter{\hbox{\epsfbox{Pics.9}}}$

\rightskip145pt

\vskip-110pt

\noindent
sense in terms of the sheaf $\Cal F|_M$. However, it does make sense
for small $\varepsilon$ because the function $f\in\Cal F|_M$ is locally
a restriction of some $\hat f\in\Cal F$. At the first glance, it may
seem that we can define $v_pf:=v_p\hat f$ even for a vector $v$ that is
not tangent to $M$ at $p\in M$. But this will not work because the
result $v_p\hat f$ will depend on the extension $\hat f$ of $f$.
The~independence of the choice of $\hat f$ is exactly the tangency of
$v$ to $M$ at a smooth point $p\in M$. Thus, we arrive at the following
intrinsic definition.

\smallskip

{\bf3.6.1.~Tangent vectors.} Let $M$ be a space with a sheaf $\Cal F$
of\break

\vskip-12pt

\rightskip0pt

\noindent
$\Bbb K$-valued functions and let $p\in M$. A $\Bbb K$-linear
functional $t:\Cal F_p\to\Bbb K$ is a {\it tangent vector\/} to $M$ at
$p$ (in~symbols, $t\in\T_pM$) if $t$ is a derivation, i.e., if
$$t(g_1g_2)=g_1(p)tg_2+g_2(p)tg_1$$
for all $g_1,g_2\in\Cal F_p$.

Let $p\in U\opensub M$. Then $(\Cal F|_U)_p=\Cal F_p$. Therefore,
assuming the induced structure on $U$, we obtain the identification
$\T_pU=\T_pM$.

For $p\in U\opensub M$, $f\in\Cal F(U)$, and $t\in\T_pM$, we define
$tf:=tf_p$.

\smallskip

{\bf3.6.2.~Exercise.} Let $t\in\T_pM$. Show that $tc=0$ for every
constant $c\in\Bbb K\subset\Cal F_p$ and that $t(\frak m_p^2)=0$.
Hence, $t$ defines a $\Bbb K$-linear functional
$\overline t:\frak m_p/\frak m_p^2\to\Bbb K$. Moreover, the
$\Bbb K$-linear map $\T_pM\to(\frak m_p/\frak m_p^2)^*$,
$t\mapsto\overline t$, is an isomorphism. By definition, $\T_pM$ and
$\T_p^*M:=\frak m_p/\frak m_p^2$ are the $\Bbb K$-linear spaces {\it
tangent\/} and {\it cotangent\/} to $M$ at $p$.

\smallskip

{\bf3.6.3.~Differential.} Let
$(M,\Cal F)\overset\psi\to\longrightarrow(N,\Cal G)$ be a smooth map
and let $p\in M$. We have the homomorphism
$\Cal F_p\overset\psi_p^*\to\longleftarrow\Cal G_{\psi(p)}$
of $\Bbb K$-algebras that is induced by the composition with $\psi$.
Hence, we get the $\Bbb K$-linear map
$\text d\psi_p:\T_pM\to\T_{\psi(p)}N$ called the {\it differential\/}
of $\psi$ at $p$. At the level of functions, the~differential is
defined via composition with $\psi$, i.e.,
$\text d\psi_pt(f):=t(f\circ\psi)$ for $f\in\Cal G(U)$,
$\psi(p)\in U\opensub N$, and $t\in\T_pM$.

We denote by
$\T M:=\bigsqcup\limits_{p\in M}\T_pM\overset\pi\to\longrightarrow M$
the disjoint union (endowed with the obvious projection) of~all tangent
spaces to points in $M$. We call $\pi:\T M\to M$ the {\it tangent
bundle\/} of $M$. Note that the {\it fibre\/} $\T_pM$ is nothing but
$\pi^{-1}(p)$.

{\unitlength=1bp$$\latexpic(40,37)(-180,0)
\put(0,31){$\T M$}\put(20,34){\vector(1,0){22}}
\put(24,25){$\text d\psi$}\put(44,31){$\T N$}
\put(10,27){\vector(0,-1){18}}\put(13,15){$\pi_M$}
\put(54,27){\vector(0,-1){18}}\put(37,15){$\pi_N$}\put(4,0){$M$}
\put(16,2){\vector(1,0){32}}\put(29,6){$\psi$}\put(50,0){$N$}
\endlatexpic$$}

\rightskip75pt

\vskip-73pt

Given a smooth map
$(M,\Cal F)\overset\psi\to\longrightarrow(N,\Cal G)$, we get the
following commutative diagram, where the differential
$\text d\psi:\T M\to\T N$ equals $\text d\psi_p$ on the fibre $\T_pM$.

\smallskip

{\bf3.6.4.~Exercise.} Show that $\T$ and $\text d$ provide a {\it
functor,} i.e., prove the following {\it chain rule.} Given smooth maps
$(L,\Cal E)\overset\varphi\to\longrightarrow(M,\Cal
F)\overset\psi\to\longrightarrow(N,\Cal G)$,
the differential of the composi-\break

\vskip-12pt

\rightskip0pt

\noindent
tion is the composition of the differentials:
$\text d(\psi\circ\varphi)=(\text d\psi)\circ(\text d\varphi)$. (The
fact that $\text d1_M=1_{\T M}$ looks quite obvious.)

\smallskip

We can picture the tangent space $\T_pM$ as the best first order
approximation of an infinitesimal neighbourhood of $p\in M$ by a
$\Bbb K$-linear space. So, the differential $\text d\psi_p$ is the best
first order linear approximation of $\psi$ over such neighbourhood.

\smallskip

{\bf3.6.5.~Tangent bundle of a subspace.} Let $S\subset M$ be a
subspace, i.e., a subset equipped with the induced structure. We denote
by $\I S$ all functions that vanish on $S$. In detail,
$\I S(U):=\big\{f\in\Cal F(U)\mid f(S\cap U)=0\big\}$ for every
$U\opensub M$. We obtain the sheaf of ideals $\I S\triangleleft\Cal F$
in the sense of the following definition.

Suppose that, for every $U\opensub M$, we are given an ideal
$\Cal J(U)\triangleleft\Cal F(U)$. We say that
$\Cal J:=\bigsqcup\limits_{U\opensub M}\Cal J(U)$ is a {\it sheaf\/} of
{\it ideals\/} in $\Cal F$ and write $\Cal J\triangleleft\Cal F$ when
the following conditions hold.

\smallskip

\leftskip60pt\rightskip60pt

\noindent
$\bullet$ If $W\opensub U\opensub M$ and $f\in\Cal J(U)$, then
$f|_W\in\Cal J(W)$.

\noindent
$\bullet$ Let us be given open subsets $U_i\opensub M$, $i\in I$, and a
function $U\overset f\to\longrightarrow\Bbb K$, where
$U:=\bigcup\limits_{i\in I}U_i$. If $f|_{U_i}\in\Cal J(U_i)$ for every
$i\in I$, then $f\in\Cal J(U)$.

\leftskip0pt\rightskip0pt

\smallskip

The germs at $p\in M$ of functions from $\Cal J$ form the {\it stalk\/}
$\Cal J_p$ of $\Cal J$ at $p$. Clearly,
$\Cal J_p\triangleleft\Cal F_p$.

\smallskip

{\bf3.6.6.~Exercise.} Let $p\in S\subset M$. Then
$(\I S)_p\subset\frak m_p$ and $(\Cal F|_S)_p=\Cal F_p/(\I S)_p$.

\smallskip

{\bf3.6.7.~Exercise.} Let $p\in S\subset M$. Show that
$\T_pS=\big\{t\in\T_pM\mid t(\I S)_p=0\big\}\le\T_pM$. This means that
the differential of the inclusion $i:S\hookrightarrow M$ can be
interpreted as an inclusion $\text di:\T S\hookrightarrow\T M$.

\smallskip

{\bf3.6.8.~Equations.} Let $(M,\Cal F)$ be a space with a sheaf of
$\Bbb K$-valued functions. One may define a closed subspace
$S\subset M$ by means of equations. Say, we could take
$E\subset\Cal F(M)$ and put
$S:=\big\{p\in M\mid e(p)=0\text{ for all }e\in E\big\}$.
Unfortunately, there are many nice spaces with sheaves where such a
definition produces nothing interesting.\footnote{Although, the
definition works somehow for the sheaves $C^\infty$.}
The reason is simple --- it can happen that $\Cal F(M)=\Bbb K$. Let us
try local functions in the equations:

Let $E\subset\Cal F$ and denote by $U_e\opensub M$ the domain of
$e\in E$, $e\in\Cal F(U_e)$. We define the subspace
$$\Z E:=\big\{p\in M\mid e(p)=0\text{ for all }e\in E\text{ such that
}p\in U_e\big\}$$
given by the {\it equations\/} $E=0$ and equipped with the induced
structure. Note that, according to this definition, $p\in\Z E$ if
$p\notin U_e$ for all $e\in E$. In particular, every closed subset is
given by equations. Indeed, let $U\opensub M$ and let
$1_U\in\Bbb K\subset\Cal F(U)$ denotes the constant $1$. Then
$M\setminus U=\Z1_U$.

We have
$\Z E=\{p\in M\mid e_p\in\frak m_p\text{ for all }e\in E\text{ such
that }p\in U_e\}$.
In particular, $\Z\Cal J=\{p\in M\mid\Cal J_p\subset\frak m_p\}$ for
any sheaf of ideals $\Cal J\triangleleft\Cal F$.

\smallskip

{\bf3.6.9.~Exercise.} Let $S\subset M$ and $E\subset\Cal F$. Show that
the operators $\Z$ and $\I$ revert the inclusion. Verify that
$\Z\I S\supset S$ and $\I\Z E\supset E$. The sheaf of ideals $\I\Z E$
is the {\it saturation\/} of $E\subset\Cal F$. Prove that the
saturation $\I\Z E$ defines the same subspace as $E$ does, i.e., that
$\Z\I\Z=\Z$. Show that $\I S$ is {\it saturated\/}, i.e., that
$\I\Z\I=\I$.

\smallskip

In order to show that (conversely) any set given by equations is
closed, we may require that the sheaf $\Cal F$ is local. A sheaf
$\Cal F$ on $M$ is {\it local\/} if every
$g\in\Cal F_p\setminus\frak m_p$ is {\it invertible\/} in $\Cal F_p$
for all $p\in M$. This means that there is some $g'\in\Cal F_p$ such
that $gg'=1$.

For every local sheaf $\Cal F$ and any $E\subset\Cal F$, the set $\Z E$
is closed in $M$. Indeed, since $\Z E=\bigcap\limits_{e\in E}\Z e$,
it~suffices to show that $\Z e$ is closed in $M$. Let $e\in\Cal F(U)$.
Then $\Z e=(M\setminus U)\cup\{p\in U\mid e_p\in\frak m_p\}$.
It~remains to prove that $\{p\in U\mid e_p\notin\frak m_p\}\opensub U$.
Let $p\in U$ and $e_p\notin\frak m_p$. Being $\Cal F$ local, we~have
$e_pf_p=1$ for suitable $p\in V\opensub M$ and $f\in\Cal F(V)$. By the
definition of germs, there exists some $W\opensub U\cap V$ such that
$p\in W$ and $e|_Wf|_W=1$. Hence, $e_qf_q=1$ for every $q\in W$. In
other words, $W\subset\{q\in U\mid e_q\notin\frak m_q\}$.

Moreover, the above arguments show that the function
$\frac1e:(M\setminus\Z e)\to\Bbb K$ defined by the rule
$p\mapsto\frac1{e(p)}$ belongs locally to $\Cal F$. So,
$\frac1e\in\Cal F(M\setminus\Z e)$ for every $e\in\Cal F$. We arrive at
another definition of a local sheaf: a sheaf $\Cal F$ is local iff, for
every $e\in\Cal F$, the locus where $e$ does not vanish is open and the
corresponding function $\frac1e$ defined on this locus belongs to
$\Cal F$.

In an arbitrary sheaf, we can sum and multiply a couple of functions
(over a locus where both are defined). In a local sheaf, we can also
perform division. Hence, it makes sense to learn how to differentiate
a fraction; by the Leibniz rule, $t\frac1g=-\frac{tg}{g^2(p)}$ for all
$t\in\T_pM$ and $g\in\Cal F_p\setminus\frak m_p$.

By Exercise 3.6.6, the sheaf $\Cal F|_S$ is local for every subspace
$S\subset M$ if $\Cal F$ is local.

\smallskip

{\bf3.6.10.~Taylor sheaves.} Suppose that every finite-dimensional
$\Bbb K$-linear space $V$ is equipped with a topology and a local sheaf
$\Cal F^V$ of $\Bbb K$-valued functions such that the following
conditions are satisfied.

\smallskip

\leftskip60pt\rightskip60pt

\noindent
$\bullet$ The topology on $V$ is the weakest one such that all
$\Cal F^V(V)\ni f:V\to\Bbb K$ are continuous.

\noindent
$\bullet$ $V^*\subset\Cal F^V(V)$.

\noindent
$\bullet$ Let $V,W$ be finite-dimensional $\Bbb K$-linear spaces. A map
$V\opensup U\overset\psi\to\longrightarrow W$ is smooth iff
$W^*\circ\psi\subset\Cal F^V(U)$.

\noindent
$\bullet$ The composition $V^*\to\frak m_p\to\frak m_p/\frak m_p^2$ is
a $\Bbb K$-linear isomorphism for every $p\in V$, where the map
$V^*\to\frak m_p$ is given by the rule
$\varphi\mapsto\varphi_p-\varphi p\in\frak m_p$.

\leftskip0pt\rightskip0pt

\smallskip

The last condition provides the identification
$\T_pV\overset\sim\to\longrightarrow V^{**}\simeq V$ given by the rule
$t\mapsto(v^*\mapsto tv^*)$, where $v^*\in V^*$. It is nothing but
Taylor's formula! Indeed, let $p\in U\opensub V$ and let
$f\in\Cal F^V(U)$. Then $f_p-f(p)\in\frak m_p$. So, there exist a
unique $\varphi\in V^*$ and $h\in\frak m_p^2$ such that
$f_p=f(p)+\varphi_p-\varphi p+h$.

In Taylor's formula, $\varphi$ provides the best linear approximation
of $f$ at $p$ modulo a term of order $2$. Hence, it is no surprise that
$\text df_pv=\varphi v$ in terms of the above identification. Indeed,
the vector $v\in V$ corresponds to the tangent vector $t\in\T_pU=\T_pV$
such that $v^*v=tv^*$ for all $v^*\in V^*$. By definition,
$\text df_pt:g\mapsto t(g\circ f)$ for all $g\in\Cal F^\Bbb K(W)$ such
that $f(p)\in W\opensub\Bbb K$. Consequently,
$\text df_pv=\text df_pt\in\T_{f(p)}\Bbb K$ corresponds to $k\in\Bbb K$
such that $k^*k=t(k^*\circ f)$ for all $k^*\in\Bbb K^*=\Bbb K$. Since
$k^*\circ f=k^*f$, we obtain
$k^*k=k^*t(f_p)=k^*t\big(f(p)+\varphi_p-\varphi
p+h\big)=k^*t\varphi_p=k^*\varphi v$,
implying $k=\varphi v$.

Let $V$ be a finite-dimensional $\Bbb K$-linear space. Then the
projection $\pi:V\oplus V\to V$ is smooth by the third and second
conditions. In particular, $U\times V=\pi^{-1}(U)\opensub(V\oplus V)$
for every $U\opensub V$. Finally, we~require that the differential
$\text df_p$ of a function depends smoothly on $p$ :

\smallskip

\leftskip60pt\rightskip60pt

\noindent
$\bullet$ Let $U\opensub V$ and let $f\in\Cal F^V(U)$. Then the
function $\text d'f:U\times V\to\Bbb K$ given by the rule
$\text d'f:(p,v)\mapsto\text df_pv$ belongs to
$\Cal F^{V\oplus V}(U\times V)$.

\leftskip0pt\rightskip0pt

\smallskip

\noindent
Sheaves $\Cal F^V$ satisfying these five conditions are called {\it
Taylor\/} sheaves.

\smallskip

There are several Taylor sheaves dealt with in geometry. The smallest
ones are formed by algebraic functions (and assume the Zariski
topology; such a topology is provided by the finite topology on
$\Bbb K$, i.e., the~weakest one with closed points). Another example is
the sheaves of analytic functions.

Here, we are interested mostly in the large sheaves $C^\infty$ of
smooth functions. Since $[v]_V\varphi$ is a constant (equal to
$\varphi v$) for any $\varphi\in V^*$, we obtain the second condition
for the sheaves $C^\infty$. Exercises 3.3.4, 3.3.5,~3.3.3, and the
solution of Exercise 3.3.3 suggested in Hints imply respectively the
first, third, fourth, and fifth conditions. It is worthwhile mentioning
that the first three conditions are valid for the sheaves $C^k$,
$k\ge0$.

\smallskip

{\bf3.6.11.~Prevarieties.} This is a crucial subsection in this
section. We want to introduce a convenient class $\widehat{\Cal V}$ of
spaces with sheaves, mostly by means of certain local properties. In
other words, every space in $\widehat{\Cal V}$ is a gluing of some
basic spaces called {\it models.} The models come from
finite-dimensional $\Bbb K$-linear spaces equipped with certain
structures.

Given Taylor sheaves $\Cal F^V$, a space $M$ with a sheaf of
$\Bbb K$-valued functions is called a {\it prevariety\/} if, locally,
it is a locally closed\footnote{A subspace $S$ in a topological space
$M$ is {\it locally closed\/} if $S=U\cap X$, where $X$ is closed in
$M$ and $U\opensub M$.}
subspace in a finite-dimensional $\Bbb K$-linear space. Usually, the
topology chosen on finite-dimensional linear spaces has a countable
basis. In order to keep this property for prevarieties, one allows only
countable or finite gluings of models (in the algebraic case, always
finite).

We denote by $\widehat{\Cal V}$ the class of all prevarieties. The
sheaves on prevarieties are obviously local. It~follows directly from
the above definition that $\widehat{\Cal V}$ is closed with respect to
taking locally closed subspaces --- called {\it subprevarieties\/} ---
and (countable or finite) gluings. A closed/open subspace in a
prevariety is called a {\it closed\/}/{\it open subprevariety.} The
intersection of finitely many (closed/open) subprevarieties is~a
(closed/open) subprevariety. Let
$(M,\Cal F)\overset\psi\to\longrightarrow(N,\Cal G)$ be a smooth map
between prevarieties and let $S\subset N$ be a (closed/open)
subprevariety. Then $\psi^{-1}(S)$ is a (closed/open) subprevariety in
$M$.

\smallskip

{\bf3.6.12.~Exercise.} Let $M\in\widehat{\Cal V}$ be a prevariety
and let $U\opensub M$. Prove that $\Cal F^M(U)$ consists of all
smooth maps $U\to\Bbb K$.

\smallskip

{\bf3.6.13.~Lemma.} {\sl For all\/ $M,N\in\widehat{\Cal V}$, there
exists a\/ $\widehat{\Cal V}$-product structure on\/ $M\times N$.}

\smallskip

{\bf Proof.} By Exercises 3.5.4 and 3.5.3, it suffices to show that
there exists a $\widehat{\Cal V}$-product structure on $V_1\times V_2$,
where the $V_i$'s are finite-dimensional linear spaces. The projection
$V_1\oplus V_2\to V_i$ is smooth by the third condition in 3.6.10. Let
$M\in\widehat{\Cal V}$ and let $\psi_i:M\to V_i$ be smooth for $i=1,2$.
We need to show that the corresponding map $\psi:M\to V_1\oplus V_2$ is
smooth. By Exercise 3.4.1, we can assume that $M$ is a model, i.e.,
$M\subset U\opensub V$, where $V$ is a finite-dimensional linear space.

Let $v_{ij}^*\in V_i^*$ be a linear basis in $V_i^*$, $i=1,2$. Then
$f_{ij}:=v_{ij}^*\circ\psi_i\in\Cal F^M(M)$ by the second condition in
3.6.10. Every function from $\Cal F^M(M)$ is locally a restriction
of a function from $\Cal F^U$. Without loss of generality, we can
therefore assume (using again Exercise 3.4.1) that
$f_{ij}=\hat f_{ij}|_M$ for all $i,j$, where
$\hat f_{ij}\in\Cal F^U(U)$. There exists a unique map
$\hat\psi_i:U\to V_i$ such that $v_{ij}^*\circ\hat\psi_i=\hat f_{ij}$
for all $j$. By the third condition in 3.6.10, $\hat\psi_i$ is smooth.
Obviously, $\psi_i=\hat\psi_i|_M$. So, we reduced the task to the case
of $M=U$. In this case, the desired fact follows immediately from the
second and third conditions in 3.6.10
$_\blacksquare$

\smallskip

We denote by $\Delta_B:=\big\{(p,p)\mid p\in B\big\}\subset B\times B$
the {\it diagonal\/} in $B\times B$. (Actually, $\Delta_B=B\times_BB$
with respect to the identity maps
$B\overset1_B\to\longrightarrow B\overset1_B\to\longleftarrow B$.)

\smallskip

{\bf3.6.14.~Lemma.} {\sl Let\/ $M_1,M_2,B\in\widehat{\Cal V}$ and let\/
$M_1\overset\varphi_1\to\longrightarrow
B\overset\varphi_2\to\longleftarrow M_2$
be smooth maps. Then the diagonal\/ $\Delta_B$ is locally closed in\/
$B\times B$. If\/ $B\subset V$ is a model, i.e., a subprevariety in a
finite-dimensional\/ $\Bbb K$-linear space\/ $V$, then $\Delta_B$ is
closed in\/ $B\times B$. There exists a fibre product\/
$\widehat{\Cal V}$-structure on\/ $M_1\times_BM_2$.}

\smallskip

{\bf Proof.} The second statement follows from
$\Delta_B=\Delta_V\cap(B\times B)$ and from
$\Delta_V=\Z_{V\times V}\{v^*\circ\pi_1-v^*\circ\pi_2\mid v^*\in
V^*\}$,
where $\pi_i:V\times V\to V$ stand for the projections.

For the first statement, we observe that
$\Delta_B\cap(B_i\times B_i)=\Delta_{B_i}$ is closed in
$B_i\times B_i$ by the second statement, where
$B=\bigcup\limits_{i\in I}B_i$ is a gluing of models $B_i\opensub B$,
$i\in I$. Therefore, $\Delta_B$ is closed in
$\bigcup\limits_{i\in I}(B_i\times\nomathbreak B_i)\opensub B\times B$.

{\unitlength=1bp$$\latexpic(40,40)(203,0)
\put(0,29){$M_1$}\put(24,32){\vector(-1,0){10}}
\put(25,29){$M_1\times M_2$}\put(66,32){\vector(1,0){13}}
\put(80,29){$M_2$}\put(7,26){\vector(0,-1){17}}\put(9,16){$\psi_1$}
\put(45,26){\vector(0,-1){17}}\put(47,16){$\psi_1\times\psi_2$}
\put(86,26){\vector(0,-1){17}}\put(88,16){$\psi_2$}\put(3,0){$B$}
\put(30,3){\vector(-1,0){17}}\put(31,0){$B\times B$}
\put(60,3){\vector(1,0){21}}\put(82,0){$B$}
\endlatexpic$$}

\leftskip108pt

\vskip-71pt

For the third statement, by Lemma 3.6.13 and Exercise 3.5.7, it
suffices to show that $M_1\times_BM_2$ is locally closed in
$M_1\times M_2$. Since $\Delta_B$ is locally closed in $B\times B$ by
the first statement, it remains to observe that
$M_1\times_BM_2=(\psi_1\times\psi_2)^{-1}(\Delta_B)$, where the map
$\psi_1\times\psi_2:M_1\times M_2\to B\times B$ in the commutative
diagram is smooth by the properties of the $\widehat{\Cal V}$-product
$B\times B$
$_\blacksquare$

\leftskip0pt

\smallskip

We are going to prove that the differential is a smooth map. First, we
need to introduce a smooth structure on the tangent bundle.

\smallskip

Let $M$ be a model. So, $M\subset U$ is a closed subprevariety in an
open subprevariety $U\opensub V$ in a finite-dimensional
$\Bbb K$-linear space $V$. We have the canonical projection
$\pi_U:\T U\to U$. The isomorphisms $\T_pU=\T_pV\simeq V$, $p\in U$,
provide the other projection $\pi':\T U\to V^{**}\simeq V$ given by
the rule\break

{\unitlength=1bp$$\latexpic(30,35)(-170,0)
\put(0,41){$\T U$}\put(17,44){\vector(1,0){20}}\put(22,37){$\sim$}
\put(39,41){$U\times V$}\put(10,37){\vector(1,-2){14}}
\put(5,18){$\pi_U$}
\put(48,37){\vector(-1,-2){14}}\put(27,0){$U$}
\endlatexpic$$}

\rightskip80pt

\vskip-78pt

\noindent
$t\mapsto(v^*\mapsto tv^*)$, where $v^*\in V^*$. Using the projections
$\pi_U,\pi'$, we get an identification $\T U\simeq U\times V$, i.e., a
{\it trivialization\/} of the tangent bundle over $U$. At the level of
fibres, this identification is an isomorphism of $\Bbb K$-linear
spaces. Since $U\times V\opensub V\oplus V$ is an open subprevariety,
we obtain the induced structure on $\T M\subset\T U\simeq U\times V$
and a smooth projection $\pi_M:\T M\to M$. By Exercise 3.6.7 and the
fifth condition in~3.6.10,

\rightskip0pt

\vskip-5pt

$$\T M={\Z}_{U\times V}(\I M\circ\pi_U)\cap{\Z}_{U\times V}\big\{\text
d'f\in\Cal F^{V\oplus V}(W\times V)\mid f\in\I M(W),\ W\opensub
U\big\}$$

{\unitlength=1bp$$\latexpic(40,36)(-187,0)
\put(5,30){$\T M$}\put(25,33){\vector(1,0){10}}
\put(37,30){$\T U$}
\put(15,27){\vector(0,-1){18}}\put(0,17){$\pi_M$}
\put(45,27){\vector(0,-1){18}}\put(32,17){$\pi_U$}\put(9,0){$M$}
\put(21,3){\vector(1,0){18}}\put(42,0){$U$}
\endlatexpic$$}

\rightskip65pt

\vskip-69pt

\noindent
is given by equations; hence, $\T M$ is closed in $\T U$ and all maps
in the commutative diagram are smooth. In other words, the structure on
$\T M$ is induced from $\T V=V\times V$ with respect to the imbedding
$M\hookrightarrow V$.

\smallskip

{\bf3.6.15.~Lemma.} {\sl Let\/ $M_i\subset U_i$ be a closed
subprevariety, where\/ $U_i\opensub V_i$ is open in a}\break

\vskip-11pt

\rightskip0pt

\noindent
{\sl finite-dimensional linear space\/ $V_i$, and let\/ $\T M_i$ be
equipped with the structure induced from\/ $\T V_i=V_i\times V_i$,
$i=1,2$. Then, for every smooth map\/ $\psi:M_1\to M_2$, the
differential\/ $\text{\rm d}\psi:\T M_1\to\T M_2$ is smooth.}

\smallskip

{\bf Proof.} We can assume that $M_2=V_2$. Let $v_j^*\in V_2^*$ be a
linear basis. The functions $v_j^*\circ\psi\in\Cal F^{M_1}(M_1)$ are
locally restrictions of some functions $f_j\in\Cal F^{U_1}$. By
Exercise 3.4.1, we can assume that $f_j\in\Cal F^{U_1}(U_1)$. There
exists a unique map $\hat\psi:U_1\to V_2$ such that
$v_j^*\circ\hat\psi=f_j$ for all $j$. In other words,
$\psi=\hat\psi|_{M_1}$. By~the third condition in 3.6.10, $\hat\psi$ is
smooth. So, we can take $M_1=U_1$.

{\unitlength=1bp$$\latexpic(40,38)(-138,0)
\put(6,30){$\T U_1$}\put(26,33){\vector(1,0){18}}
\put(28,36){$\text d\psi$}\put(46,30){$V_2\times V_2$}
\put(80,33){\vector(1,0){12}}\put(82,35){$\pi'$}\put(94,30){$V_2$}
\put(15,27){\vector(0,-1){18}}\put(0,17){$\pi_{U_1}$}
\put(62,27){\vector(0,-1){18}}\put(47,17){$\pi_{V_2}$}
\put(98,27){\vector(0,-1){18}}\put(88,17){$v^*$}\put(11,0){$U_1$}
\put(22,3){\vector(1,0){34}}\put(36,6){$\psi$}\put(58,0){$V_2$}
\put(94,0){$\Bbb K$}
\endlatexpic$$}

\rightskip118pt

\vskip-68pt

By the properties of the $\widehat{\Cal V}$-product $V_2\times V_2$, it
suffices to show that $\pi'\circ\text d\psi:\T U_1\to V_2$ is smooth
because $\pi_{V_2}\circ\text d\psi=\psi\circ\pi_{U_1}$ is smooth. By
the third condition in 3.6.10, we need only to verify that
$v^*\circ\pi'\circ\text d\psi\in\Cal F^{U_1\times V_1}(U_1\times V_1)$
for every $v^*\in V_2^*$. Hence, by the fifth condition in 3.6.10, it
remains to check that $v^*\circ\pi'\circ\text d\psi=\text d'f$, where
$f:=v^*\circ\psi\in\Cal F^{U_1}(U_1)$.

\rightskip0pt

Let $p\in U_1$, let $t\in\T_pU_1$, and let $v\in V_1$, $v'\in V_2$
be the vectors corresponding to $t$, $\text d\psi t$. This means that
$t\varphi=\varphi v$ for all $\varphi\in V_1^*$, that
$\pi'(\text d\psi t)=v'$, and that $(\text d\psi t)v^*=v^*v'$. By the
fourth condition in~3.6.10, we have $f_p=f(p)+\varphi_p-\varphi(p)+h$
with $h\in\frak m_p^2\subset\Cal F_p^{U_1}$ and $\varphi\in V_1^*$.
Consequently,
$$\text d'f(p,v)=\text df_pv=\varphi
v=t\varphi=tf=t(v^*\circ\psi)=(\text d\psi t)v^*=v^*v'=
v^*\big(\pi'(\text d\psi t)\big)=(v^*\circ\pi'\circ\text d\psi)t\
_\blacksquare$$

\smallskip

Taking $M_1:=M_2:=M$ and $\psi:=1_M$ in Lemma 3.6.15, we can see that
the induced structure on $\T M\subset\T V$ is independent of the choice
of an embedding $M\hookrightarrow V$ into a linear space.

\smallskip

Let $M$ be an arbitrary prevariety. It is a gluing of models
$M=\bigcup\limits_{i\in I}M_i$. By Exercise 3.4.2, we can introduce a
structure on $\T M$ as a gluing of the structures on
$\T M_i\subset\T M$ because the structures on $\T(M_i\cap M_j)$ induced
from $\T M_i$ and from $\T M_j$ are the same by Lemma 3.6.15. A similar
argument shows that the structure constructed on $\T M$ is independent
of the choice of a gluing $M=\bigcup\limits_{i\in I}M_i$.
By~Exercise~3.4.1 and Lemma 3.6.15, the differential $\text d\psi$ of a
smooth map $\psi:M_1\to M_2$ between prevarieties is a smooth map.

\smallskip

{\bf3.6.16.~Exercise.} Let $M\in\widehat{\Cal V}$ be a prevariety.
Show that the maps $\T M\times_M\T M\overset+\to\longrightarrow\T M$,
$(t_1,t_2)\mapsto t_1+t_2$, and
$\Bbb K\times\T M\overset\cdot\to\longrightarrow\T M$,
$(k,t)\mapsto kt$, are smooth. In words, the operations $+$ and $\cdot$
are smooth on the tangent bundle (where defined).

\medskip

{\bf3.7.~$C^\infty$-manifolds.} Let $M$ be a hausdorff topological
space equipped with a sheaf of $\Bbb K$-valued $C^\infty$-functions and
possessing a countable basis of topology. We say that $M$ is a
$C^\infty$-{\it manifold\/} (or simply a {\it manifold\/}) if, locally,
it is an open subvariety in a finite-dimensional $\Bbb K$-linear space.

Let $T_1,T_2$ be topological spaces. The weakest topology on
$T_1\times T_2$ with continuous projections
$\pi_i:T_1\times T_2\to T_i$ is called the {\it product\/} topology. We
must warn the reader that the topology introduced in Subsections 3.5
and 3.6.11 on $\widehat{\Cal V}$-products may be stronger than the
product topology as it happens, for instance, in the case of the
sheaves of algebraic functions. However, for the sheaves $C^\infty$,
these topologies coincide by Exercise 3.3.4.

\smallskip

{\bf3.7.1.~Exercise.} Show that a topological space $T$ is hausdorff
iff the diagonal $\Delta_T$ is closed in the space $T\times T$ equipped
with the product topology.

\smallskip

{\bf3.7.2.~Families and bundles.} Let $\pi_i:T_i\to B$ be smooth maps
between prevarieties $T_i,B\in\widehat{\Cal V}$, $i=1,2$. We can
interpret $\pi_i$ as a {\it family\/} of spaces $\pi_i^{-1}(p)$, called
{\it fibres,} parameterized by $p\in B$. A~{\it morphism\/} between
such families is a smooth map $\psi:T_1\to T_2$ such that
$\pi_2\circ\psi=\pi_1$. Obviously, the~composition of morphisms is a
morphism and the identity map is a morphism. An invertible morphism (=
possessing a two-side inverse) is an {\it isomorphism.}

Let $F,B\in\widehat{\Cal V}$ be prevarieties. A {\it trivial\/}
({\it fibre\/}) {\it bundle\/} over $B$ is a family of subspaces
$\pi:T\to B$ isomorphic to the trivial family $F\times B\to B$. In
other words, a trivial bundle is a product that has forgotten one of
its projections. A family of subspaces $\pi:T\to B$ is a
({\it fibre\/}) {\it bundle\/} if it is {\it locally trivial,} i.e., if
there exists an open cover of the {\it base\/}
$B=\bigcup\limits_{i\in I}B_i$, called a {\it trivializing\/} cover,
such that $\pi^{-1}(B_i)\to B_i$ is a trivial bundle for all $i\in I$.
It is immediate that a bundle over a manifold whose fibres are
manifolds is a manifold. As we have seen in Subsection 3.6.11, the
tangent bundle $\pi:\T M\to M$ of any manifold $M$ is a bundle.
However, in general, the tangent bundle of a prevariety is not a
bundle!\footnote{For the sheaves $C^\infty$, the tangent bundle of a
prevariety which is not a manifold can be a bundle (take, for example,
a closed ball). In the case of algebraic geometry, the tangent bundle
of a prevariety is rarely a bundle. This happens, say, when the
prevariety is smooth and rational.}
A bundle with discrete fibres is called a ({\it regular\/}) {\it
covering.} Coverings are essential when studying manifolds that carry a
geometrical structure (see Section 5). The reader can see the picture
of a simple covering at the very beginning of Section 2.

\smallskip

{\bf3.7.3.~Exercise.} Prove that the sphere and the projective space
are compact manifolds.

\smallskip

{\bf3.7.4A.~Example.} More generally, prove that the grassmannians
$\Gr_\Bbb K(k,V)$ and $\Gr_\Bbb R^+(k,V)$ are compact manifolds (see
2.12A and 3.4.4A).

\smallskip

{\bf3.7.5.~Exercise.} Let $V$ be a finite-dimensional $\Bbb K$-linear
space. Show that
$$T:=\big\{(l,v)\mid V\ge l\ni v,\ \dim_\Bbb Kl=1\big\}\subset\Bbb
P_\Bbb KV\times V$$
is a closed submanifold and that the projection to $\Bbb P_\Bbb KV$
provides a bundle $\pi:T\to\Bbb P_\Bbb KV$. This bundle is called {\it
tautological.} Visualize $T$ as a M\"obius band in the case of
$\dim_\Bbb RV=1$. Is every tautological bundle trivial?

\smallskip

{\bf3.7.6A.~Exercise.} More generally, formulate and solve a similar
exercise about grassmannians.

\smallskip

{\bf3.7.7.~Exercise.} Prove that the surface of a $3$-cube in
$\Bbb R^3$ is not a $C^\infty$-manifold.

\smallskip

{\bf3.7.8.~Tangent vector to a curve.} A smooth map
$\Bbb R\opensup(a,b)\overset c\to\longrightarrow M$ into a
$C^\infty$-prevariety $(M,\Cal F)$ is a {\it parameterized\/} smooth
{\it curve.} The {\it tangent vector\/} $\dot c(t_0)$ to the curve $c$
at the point $c(t_0)$ is given by the formula
$\Cal F_{c(t_0)}\ni f_{c(t_0)}\mapsto\frac{\text d}{\text
dt}\big|_{t=t_0}f\big(c(t)\big)$.
It is easy to see that every tangent vector to a manifold is tangent to
a suitable smooth curve.

\smallskip

{\bf3.7.9.~Exercise.} Translate any book on basic differential
topology (the worst is the best) into the terms of the above
exposition.

\medskip

{\bf3.8A.~Final remarks.} It is important to study not only smooth
manifolds but also manifolds with singularities (analytic spaces in the
case of analytic sheaves). Such hausdorff spaces --- let us call them
{\it varieties\/} --- should be defined by means of models
$M\subset U\opensub V$ whose sheaf $\I M$ of ideals satisfies certain
finiteness conditions. In this case, our considerations in 3.6.11--16
should work for varieties.

There are indications that a right definition of a smooth space should
be close to the one mentioned in Remark 3.8.1A below. However, if we
were to simply accept it, we would not have had the above journey
around the world of smooth spaces.

\smallskip

{\bf3.8.1A.~Remark.} Let $(M,\Cal F^M)$ and $(\T^*,\Cal F^{\T^*})$ be
spaces with sheaves of $\Bbb K$-valued functions and let
$\pi:\T^*\to M$ be a smooth map whose fibres are finite-dimensional
$\Bbb K$-linear spaces such that the global operations
$+:\T^*\times_M\T^*\to\T^*$ and $\cdot:\Bbb K\times\T^*\to\T^*$ are
smooth. It seems possible to define varieties in these terms by using a
de Rham morphism of the sheaves $\text d:\Cal F^M\to\Cal T^*$ subject
to a Leibniz rule, where $\Cal T^*$ stands for the sheaf of smooth
sections of $\pi:\T^*\to M$.

\rightheadtext{Elementary geometry}

\bigskip

\centerline{\bf4.~Elementary geometry}

\medskip

\rightline{\oitossi there were and are even now geometers and
philosophers}

\rightline{\oitossi$\dots$ who doubt that the whole
universe $\dots$ was created}

\rightline{\oitossi purely in accordance with Euclidean
geometry}

\medskip

\rightline{{\oitoss --- FYODOR DOSTOYEVSKY,} {\oitossi The Karamazov
brothers}}

\bigskip

\rightline{\oitossi Out of nothing I have created a strange new
universe.}

\medskip

\rightline{{\oitoss --- J\'ANOS BOLYAI}}

\bigskip

For a long time, there was little doubt that Euclidean geometry is the
`right' geometry; nowadays, non-Euclidean geometry is involved in many
areas of mathematics and physics. It is no exaggeration to say that the
discovery of non-Euclidean geometry, more specifically of hyperbolic
geometry, represented a major mathematical and philosophical
breakthrough. The ancient question concerning the fifth
postulate\footnote{Roughly speaking, the postulate says: given a point
and a line, there exists a unique parallel line passing through the
point.}
was finally answered, and the answer was astonishing: the apparently
evident fifth postulate turned out to be {\it independent\/} since
hyperbolic and Euclidean geometries share the same axioms except the
fifth one (which is false in the hyperbolic plane). Of course, we are
not interested in axiomatic geometry here. Instead, we study hyperbolic
and many other non-Euclidean geometries on the basis of simple linear
algebra. In this regard, the reader is welcome to consult Section 6
devoted to linear and hermitian tools.

\medskip

{\bf4.1.~Some notation.} Let $V$ be a finite-dimensional
$\Bbb K$-linear space equipped with a nondegenerate hermitian form
$\langle-,-\rangle$, where $\Bbb K=\Bbb R$ or $\Bbb K=\Bbb C$.
Depending on the context, we will frequently use a same letter to
denote a point in $\Bbb P_\Bbb KV$ and a representative in $V$. We use
the notation and convention for projectivizations introduced in
Subsection 2.6 : given a subset $S\subset V$, the image of $S$ under
the quotient map $\pi:V^\centerdot\to\Bbb P_\Bbb KV$ is denoted by
$\Bbb P_\Bbb KS:=\pi\big(S\setminus\{0\}\big)\subset\Bbb P_\Bbb KV$.

The {\it signature\/} of $p\in\Bbb P_\Bbb KV$ is the sign of
$\langle p,p\rangle$ (it can be $-$, $+$, or $0$). Note that signature
is well defined since, for another representative $kp\in V$,
$k\in\Bbb K^\centerdot$, we have
$\langle kp,kp\rangle=|k|^2\langle p,p\rangle$. The projective space
$\Bbb P_\Bbb KV$ is divided into three disjoint parts consisting of
{\it negative,} {\it positive,} and {\it isotropic\/} points:
$$\B V:=\{p\in\Bbb P_\Bbb KV\mid\langle p,p\rangle<0\},\qquad\E
V:=\{p\in\Bbb P_\Bbb KV\mid\langle p,p\rangle>0\},\qquad\S
V:=\{p\in\Bbb P_\Bbb KV\mid\langle p,p\rangle=0\}.$$
The isotropic points constitute the {\it absolute\/} $\S V$ of
$\Bbb P_\Bbb KV$. The absolute is a `wall' separating the geometries
(not yet introduced) on $\B V$ and $\E V$. Moreover, we will see later
that the absolute itself possesses its own geometry. We denote
$\overline\B V:=\B V\sqcup\S V$ and $\overline\E V:=\E V\sqcup\S V$.

Let $p\in\Bbb P_\Bbb KV$ be nonisotropic. We introduce the following
notation for the orthogonal decomposition:
$$V=\Bbb Kp\oplus p^\perp,\qquad v=\pi'[p]v+\pi[p]v,$$
where
$$\pi'[p]v:=\frac{\langle v,p\rangle}{\langle p,p\rangle}p\in\Bbb
Kp,\qquad\pi[p]v:=v-\frac{\langle v,p\rangle}{\langle p,p\rangle}p\in
p^\perp.$$
It is easy to see that $\pi'[p]$ and $\pi[p]$ do not depend on the
choice of a representative $p\in V$.

\medskip

{\bf4.2.~Tangent space.} Let $p\in\Bbb P_\Bbb KV$, let $f$ be a smooth
function defined on an open neighbourhood $U\opensub\Bbb P_\Bbb KV$ of
$p$, and let $\varphi:\Bbb Kp\to V$ be a $\Bbb K$-linear map. Using the
notation from Subsection 3.3, we~define
$$t_\varphi f:=(\varphi p)_p\tilde f,$$
where $\tilde f$ stands for the lift of $f$ to an open neighbourhood of
$\Bbb K^\centerdot p$ in $V$. This lift satisfies
$\tilde f(kp)=\tilde f(p)$ for all $k\in\Bbb K^\centerdot$.

\smallskip

{\bf4.2.1.~Exercise.} Verify that $t_\varphi$ is well defined and
conclude that $t_\varphi\in\T_p\Bbb P_\Bbb KV$. Show that $t_\varphi=0$
iff $\varphi p\in\Bbb Kp$. Therefore,
$\T_p\Bbb P_\Bbb KV=\Lin_\Bbb K(\Bbb Kp,V/\Bbb Kp)$. For a nonisotropic
$p\in\Bbb P_\Bbb KV$, we have the identifications
$\T_p\Bbb P_\Bbb KV=\Lin_\Bbb K(\Bbb Kp,p^\perp)=\langle-,p\rangle
p^\perp$,
where $\langle-,p\rangle v:x\mapsto\langle x,p\rangle v$.

\smallskip

\vskip10pt

\noindent
\hskip410pt$\vcenter{\hbox{\epsfbox{Pics.11}}}$

\rightskip50pt

\vskip-112pt

Intuitively, we can interpret the identification
$\T_p\Bbb P_\Bbb KV=\Lin_\Bbb K(\Bbb Kp,p^\perp)$ as follows. A point
$p\in\Bbb P_\Bbb KV$ corresponds to a line $L\subset V$ passing through
$0$. A tangent vector $t_\varphi$ at $p$ is an infinitesimal movement
of $L$ (a sort of rotation about $0$) and so can be exhibited as a
direction orthogonal to $L$. But this direction is not merely an
element $t_\varphi p\in p^\perp$ : the fact that $t_\varphi$ is a
linear map provides the independence of the choice of a representative
$p\in V$.

\smallskip

The tangent vector to a smooth curve in $\Bbb P_\Bbb KV$ at a
nonisotropic $p$ can be handy expressed in terms of the identification
$\T_p\Bbb P_\Bbb KV=\Lin_\Bbb K(\Bbb Kp,p^\perp)$ :

\smallskip

{\bf4.2.2.~Exercise.} Let $c:(a,b)\to\Bbb P_\Bbb KV$ be a smooth curve,
let $c_0:(a,b)\to V$ be a smooth lift of $c$ to~$V$, and let $c(t)$ be
a nonisotropic point, $t\in(a,b)$. Show that the tangent vector to $c$
at $c(t)$ corresponds to the $\Bbb K$-linear map
$\dot c(t):\Bbb Kc_0(t)\to c_0(t)^\perp$,
$c_0(t)\mapsto\pi\big[c_0(t)\big]\dot c_0(t)$.

\rightskip0pt

\smallskip

{\bf4.2.3.*~Exercise.} Let $W\le V$ be an $\Bbb R$-linear subspace. A
point $p\in W$ is said to be {\it projectively smooth\/} if
$\dim_\Bbb R(\Bbb Kp\cap W)=\min\limits_{0\ne w\in W}\dim_\Bbb R(\Bbb
Kw\cap W)$.
Prove that the projectivization $\Bbb P_\Bbb KS\subset\Bbb P_\Bbb KV$
of the subset $S\subset W$ formed by all projectively smooth points in
$W$ is a submanifold. Let $p\in S$ be a projectively smooth point and
let $\varphi:\Bbb Kp\to V$ be a $\Bbb K$-linear map. Show that
$t_\varphi\in\T_p\Bbb P_\Bbb KS$ iff $\varphi p\in W+\Bbb Kp$.

\medskip

{\bf4.3.~Metric.} Let $p\in\Bbb P_\Bbb KV$ be a nonisotropic point.
Given $v\in p^\perp$, we define
$$t_{p,v}:=\langle-,p\rangle v\in{\T}_p\Bbb P_\Bbb KV.$$
Note that $t_{p,v}$ does depend on the choice of a representative
$p\in V$ : If we pick a new representative $\overline kp\in V$,
$\overline k\in\Bbb K^\centerdot$, then we must take $\frac1kv\in V$ in
place of $v$ in order to keep $t_{p,v}$ the same.

The tangent space $\T_p\Bbb P_\Bbb KV$ is equipped with the hermitian
form
$$\langle t_{p,v_1},t_{p,v_2}\rangle:=\pm\langle p,p\rangle\langle
v_1,v_2\rangle.\leqno{\bold{(4.3.1)}}$$
This definition is correct as the formula is independent of the choice
of representatives $p,v_1,v_2\in V$ providing the same
$t_{p,v_1},t_{p,v_2}$. One can readily see that this hermitian form,
called a {\it hermitian metric\/} (or simply a {\it metric\/}), depends
smoothly on a nonisotropic $p$. Actually, this is another instance of a
typical situation when we are to show some smooth dependence on a
parameter. In general, such cases can be treated as in Exercise 3.6.16
and usually the concept of fibred product is to be explored. The~only
essential step in the proof consists in observing the (say) algebraic
nature of the formulae involving the parameter.

\smallskip

What do we need a hermitian metric for?

\smallskip

{\bf4.3.2.~Length and angle.} Let $M$ be a smooth manifold such that
every tangent space $\T_pM$ is equipped with a positive-definite
hermitian form $\langle-,-\rangle$ depending smoothly on $p$. Then we
can measure the length of a smooth curve $c:[a,b]\to M$ by using the
familiar formula
$$\ell c:=\int_a^b\sqrt{\big\langle\dot c(t),\dot
c(t)\big\rangle}\,\text dt,$$
where $\dot c(t)$ stands for the tangent vector to $c$ at $c(t)$.

We can also measure the nonoriented angle $\alpha\in[0,\pi]$ between
nonnull tangent vectors $0\ne t_1,t_2\in\T_pM$ by using the other
familiar formula
$$\cos\alpha=\frac{\Re\langle t_1,t_2\rangle}{\sqrt{\langle
t_1,t_1\rangle}\cdot\sqrt{\langle t_2,t_2\rangle}}.$$

In the particular case when $\Bbb K=\Bbb C$ and the real subspace
$\Bbb Rt_1+\Bbb Rt_2\le\T_pM$ is complex, the oriented angle
$\alpha\in[0,2\pi)$ from $t_1$ to $t_2$ is given by
$\alpha=\Arg\langle t_2,t_1\rangle$.

\smallskip

In other words, a {\it hermitian metric\/} is what equips the manifold
with a geometric structure.

\medskip

{\bf4.4.~Examples.} By taking a particular field $\Bbb K$ and a
signature of the form $\langle-,-\rangle$ on $V$, we get many examples
of classic geometries.

\smallskip

$\bullet$ We take $\Bbb K=\Bbb C$, $\langle-,-\rangle$ of signature
$++$, and the sign $+$ in (4.3.1). The Riemann sphere $\Bbb P_\Bbb CV$
becomes a {\it round\/} sphere. It looks just like the usual sphere (of
radius $\frac12$) in Euclidean $3$-dimensional space (see 4.5.4).

\smallskip

$\bullet$ We take $\Bbb K=\Bbb C$, $\langle-,-\rangle$ of signature
$-+$, and the sign $-$ in (4.3.1).

\smallskip

{\bf4.4.1.~Exercise.} Show that the Riemann sphere $\Bbb P_\Bbb CV$ is
formed by the closed discs $\overline\B V$ and $\overline\E V$ glued
along the absolute $\S V$. Note that the hermitian metric on
$\T_p\Bbb P_\Bbb CV$ is positive-definite for all nonisotropic $p$.

\smallskip

Each of $\B V$ and $\E V$ is a {\it Poincar\'e\/} disc. It is endowed
with the corresponding metric and constitutes the most famous model of
plane {\it hyperbolic\/} geometry. We call $\Bbb P_\Bbb CV$ the {\it
Riemann-Poicar\'e\/} sphere.\footnote{We thank Pedro Walmsley Frejlich
for suggesting this term.}

\smallskip

$\bullet$ We take $\Bbb K=\Bbb R$, $\langle-,-\rangle$ of signature
$-++$, and the sign $-$ in (4.3.1).

\smallskip

{\bf4.4.2.~Exercise.} Show that the real projective plane
$\Bbb P_\Bbb RV$ is formed by the closed disc $\overline\B V$ and
M\"obius band $\overline\E V$ glued along the absolute $\S V$. Note
that the metric on $\T_p\Bbb P_\Bbb RV$ is positive-definite for
$p\in\B V$ and has signature $-+$ for $p\in\E V$.

\smallskip

The metric on the M\"obius band $\E V$ is not positive-definite (it is
called a {\it lorentzian\/} metric). In~spite of this fact, the metric
still equips $\E V$ with its adequate geometry. The fact that the
concepts of length and angle do not work fairly in this case does not
mean at all that the geometry has been lost (see~Subsection 4.5.11).

The disc $\B V$ equipped with its metric is known as the {\it
Beltrami-Klein\/} disc. It constitutes another model of plane
hyperbolic geometry. It is easy to show (see Exercise 4.5.10) that the
Beltrami-Klein disc and the Poincar\'e disc are essentially isometric.
However, there is something fundamentally different about these two
hyperbolic spaces: while the complement of a Poincar\'e disc in
$\Bbb P_\Bbb CV$ is another Poincar\'e disc, the~complement of the\
Beltrami-Klein disc in $\Bbb P_\Bbb RV$ is a lorentzian M\"obius band
$\dots$ we will soon discover that there is more to the above sentence
than just naming five great mathematicians.

\smallskip

$\bullet$ We take $\Bbb K=\Bbb C$, $\langle-,-\rangle$ of signature
$-++$, and the sign $-$ in (4.3.1). The open $4$-ball
$\B V\subset\Bbb P_\Bbb CV$ is the {\it complex hyperbolic\/} plane.
We call the entire $\Bbb P_\Bbb CV$ the {\it extended complex
hyperbolic\/} plane. It is curious that all the above examples can be
naturally embedded into the extended complex hyperbolic plane
(see~4.7). Moreover, one can deform an embedded round sphere into a
Riemann-Poincar\'e sphere~$\dots$ Which geometry should appear along
the way of the deformation?

\smallskip

$\bullet$ We take $\Bbb K=\Bbb R$, $\langle-,-\rangle$ of signature
$-+++$, and the sign $-$ in (4.3.1). The open $3$-ball
$\B V\subset\Bbb P_\Bbb RV$ is the {\it real hyperbolic\/} space. The
manifold $\E V$ --- called the {\it de Sitter\/} space --- is
lorentzian, i.e., the~signature of the metric on $\T_p\E V$ is $-++$
for all $p\in\E V$. The de Sitter space is popular among physicists as
they think it applies to general relativity.

\smallskip

$\bullet$ We take $\Bbb K=\Bbb C$, $\langle-,-\rangle$ of signature
$+\dots+$, and the sign $+$ in (4.3.1). We get the projective space
$\Bbb P_\Bbb CV$ equipped with the positive-definite {\it
Fubini-Study\/} metric. This metric is essential in many areas of
mathematics and physics, including complex analysis and
classical/quantum mechanics.

\medskip

{\bf4.5.~Geodesics and tance.} Let $W\le V$ be a $2$-dimensional
$\Bbb R$-linear subspace such that the hermitian form, being
restricted to $W$, is real and nonnull. We call
$\Bbb P_\Bbb KW\subset\Bbb P_\Bbb KV$ a {\it geodesic.}

\smallskip

{\bf4.5.1.~Exercise.} Show that $\Bbb Kp\cap W=\Bbb Rp$ for all
$0\ne p\in W$ and that $\Bbb P_\Bbb KW=\Bbb P_\Bbb RW$. Hence, every
geodesic is topologically a circle. The geodesic $\Bbb P_\Bbb KW$ spans
its projective line $\Bbb P_\Bbb K(\Bbb KW)\subset\Bbb P_\Bbb KV$.
The~geodesics $\Bbb P_\Bbb KW_1$ and $\Bbb P_\Bbb KW_2$ are equal iff
$W_1=kW_2$ for some $k\in\Bbb K^\centerdot$.

\smallskip

{\bf4.5.2.~Exercise.} Let $\Bbb P_\Bbb KV$ be a projective line,
$\dim_\Bbb KV=2$. Given a nonisotropic $p\in\Bbb P_\Bbb KV$, there
exists a unique $q\in\Bbb P_\Bbb KV$ such that $\langle p,q\rangle=0$
(in words, $q$ is orthogonal to $p$). Let $p_1,p_2\in\Bbb P_\Bbb KV$ be
distinct points. If $p_1,p_2$ are nonorthogonal, then there exists a
unique geodesic containing $p_1,p_2$. If~$\langle p_1,p_2\rangle=0$ and
$p_1$ is nonisotropic, then every geodesic in $\Bbb P_\Bbb KV$ passing
through $p_1$ passes also through $p_2$.

\smallskip

{\bf4.5.3.~Exercise.} Let $p\in\Bbb P_\Bbb KV$ be a nonisotropic point
and let $0\ne t\in\T_p\Bbb P_\Bbb KV$ be a nonnull tangent vector at
$p$. Show that there exists a unique geodesic passing through $p$ with\
tangent vector $t$. Let~$p_1,p_2\in\Bbb P_\Bbb KV$ be distinct
nonorthogonal points with nonisotropic $p_1$ and let $G$ be the
geodesic that passes through $p_1$ and $p_2$. We denote by $q\in G$
the point orthogonal to $p_1$. Show that
$\langle-,p_1\rangle\frac{\pi[p_1]p_2}{\langle p_2,p_1\rangle}$ is a
tangent vector at $p_1$ to the oriented segment of geodesic from $p_1$
to $p_2$ not passing through $q$.

\smallskip

Let us calculate the length of geodesics. By Exercise 4.5.1, we can
assume that $\dim_\Bbb KV=2$.

\smallskip

{\bf4.5.4.~Spherical geodesics.} A geodesic $\Bbb P_\Bbb KW$ is {\it
spherical\/} if $W$ has signature $++$. Such a geodesic spans the
projective line $\Bbb P_\Bbb KV$ with $V$ of signature $++$. We will
parameterize $\Bbb P_\Bbb KW$. Let $p_1\in W$. We~include $p_1$ in an
orthonormal basis $p_1,q\in V$ with $q\in W$. The curve
$$c_0:[0,a]\to V,\qquad c_0(t):=p_1\cos t+q\sin t,\qquad a\ge0$$
is a lift to $V$ of a segment of geodesic $c:[0,a]\to\Bbb P_\Bbb KV$
joining $p_1=c(0)$ and $p_2:=c(a)$. By Exercise~4.2.2, the tangent
vector to $c$ at $c(t)$ equals
$$\dot c(t)=\big\langle-,c_0(t)\big\rangle\frac{\pi\big[c_0(t)\big]
\dot c_0(t)}{\big\langle
c_0(t),c_0(t)\big\rangle}=\big\langle-,c_0(t)\big\rangle\,\dot
c_0(t)$$
because $\big\langle c_0(t),c_0(t)\big\rangle=1$ and
$\big\langle\dot c_0(t),c_0(t)\big\rangle=0$ for all $t\in[0,a]$.
Hence,
$\ell c=\int\limits_0^a\sqrt{\big\langle\dot c(t),\dot
c(t)\big\rangle}\,\text dt=\int\limits_0^a\,\text dt=\nomathbreak a$
(we take the sign $+$ in (4.3.1)). If $a\in[0,\frac\pi2]$, then $a$ can
be expressed in terms of the {\it tance}
$$\ta(p_1,p_2):=\frac{\langle p_1,p_2\rangle\langle
p_2,p_1\rangle}{\langle p_1,p_1\rangle\langle p_2,p_2\rangle}.
\leqno{\bold{(4.5.5)}}$$
By Sylvester's criterion, $\ta(p_1,p_2)\in[0,1]$ with the extremal
values corresponding to $p_2=q$ and $p_2=p_1$. A direct calculation
shows that $\ta(p_1,p_2)=\cos^2a$. Therefore,
$$\ell c=\arccos\sqrt{\ta(p_1,p_2)}.$$

Let $p_1,p_2\in\Bbb P_\Bbb KW$ be distinct nonorthogonal points in a
spherical geodesic. They divide the circle $\Bbb P_\Bbb KW$ into two
segments. The one that does not contain the point orthogonal ({\it
antipodal\/}) to $p_1$ is the {\it shortest\/} segment joining $p_1$
and $p_2$ and its length $a<\frac\pi2$ is given by the above formula.
When $p_1,p_2$ are orthogonal, either segment has length $\frac\pi2$.
This is why, in Examples 4.4, the round sphere has radius $\frac12$.

\smallskip

{\bf4.5.6.~Hyperbolic geodesics.} A geodesic $\Bbb P_\Bbb KW$ is {\it
hyperbolic\/} if $W$ has signature $-+$. Such a geodesic spans the
projective line $\Bbb P_\Bbb KV$ with $V$ of signature $-+$. We will
parameterize $\Bbb P_\Bbb KW$. Let~$p_1\in W$ be nonisotropic. We
include $p_1$ in an orthonormal basis $p_1,q\in V$ with $q\in W$. The
curve
$$c_0:[0,a]\to V,\qquad c_0(t):=p_1\cosh t+q\sinh t,\qquad a\ge0$$
is a lift to $V$ of a segment of geodesic $c:[0,a]\to\Bbb P_\Bbb KV$
joining $p_1=c(0)$ and $p_2:=c(a)$. (The~hyperbolic functions are
defined as $\cosh t:=\frac{e^t+e^{-t}}2$ and
$\sinh t:=\frac{e^t-e^{-t}}2$.) It is easy to see that
$\big\langle c_0(t),c_0(t)\big\rangle=\langle p_1,p_1\rangle$ for all
$t\in[0,a]$. So, the segment $c$ contains no isotropic points. As
above, $\ell c=a$ (we take the sign $-$ in (4.3.1)). By Sylvester's
criterion, $\ta(p_1,p_2)\ge1$ with the extremal value corresponding to
$p_2=p_1$. Hence,
$$\ell c=\arccosh\sqrt{\ta(p_1,p_2)}.$$

A hyperbolic geodesic contains exactly two isotropic points called {\it
vertices.} They divide the geodesic into two parts; one is positive and
the other, negative. The vertices can be treated as points at infinity.

\smallskip

{\bf4.5.7.~Triangle inequality.} We can use the above expressions and
introduce distance functions in the parts of $\Bbb P_\Bbb KV$ where the
hermitian metric (4.3.1) is positive-definite: the {\it hyperbolic
distance\/} $d(p_1,p_2):=\arccosh\sqrt{\ta(p_1,p_2)}$ is a distance
function in the real or complex hyperbolic geometries; the~{\it
spherical distance\/} $d(p_1,p_2):=\arccos\sqrt{\ta(p_1,p_2)}$ is a
distance function in the Fubini-Study spaces.

These formulae are monotonic in tance. Therefore, it sounds like a good
idea to use tance in place of distance because tance is a simple
algebraic expression (involving just the hermitian form on $V$
which~is, after all, the source  of the geometry on $\Bbb P_\Bbb KV$).
We know that distance is additive. Better to say, it is subject to the
triangle inequality. Let us express this inequality in terms of tances.

We consider the real hyperbolic case. Take $\Bbb K=\Bbb R$,
$\langle-,-\rangle$ of signature $-+++$, and the sign $-$ in (4.3.1).
Let $p_1,p_2,p_3\in\B V$. We fix representatives such that
$\langle p_i,p_i\rangle=-1$ and $r_1,r_2>0$, where
$r_i:=-\langle p_i,p_{i+1}\rangle$ (the indices are modulo
$3$). By Sylvester's criterion, $r_i^2\ge1$ and
$\det\left[\smallmatrix-1&-r_1&-r_3\\-r_1&-1&-r_2\\-r_3&-r_2&-1
\endsmallmatrix\right]\le0$. Hence,
$$r_1^2+r_2^2+r_3^2\le2r_1r_2r_3+1,\leqno{\bold{(4.5.8)}}$$
implying $r_3\ge1$. The triangle inequality
$\arccosh r_1\le\arccosh r_2+\arccosh r_3$ is equivalent to
$$r_1\le\cosh(\arccosh r_2+\arccosh
r_3)=r_2r_3+\sqrt{r_2^2-1}\sqrt{r_3^2-1}$$
(since $\cosh(x+y)=\cosh x\cosh y+\sinh x\sinh y$) and follows from
$(r_1-r_2r_3)^2\le(r_2^2-1)(r_3^2-1)$. We arrived at (4.5.8). The
inequality (4.5.8) is the triangle inequality in terms of tances. It
codifies simultaneously the three triangle inequalities involving
$p_1,p_2,p_3$. The equality occurs exactly when $p_1,p_2,p_3$ belong to
a same geodesic.

\smallskip

{\bf4.5.9.~Exercise.} Prove the triangle inequalities for
the complex hyperbolic plane $\B V$ and for the Fubini-Study spaces.

\smallskip

In conclusion: there is no need to deal with distances in the hermitian
manifolds under consideration. All we need is tance, hermitian algebra,
and the synthetic description of geodesics introduced above. The fact
(sometimes taken as a definition) that a geodesic is a curve locally
minimizing distance is of course valid in our case. We postpone the
proof of this fact until Appendix 10A.

\smallskip

{\bf4.5.10.~Exercise.} Identify the Poincar\'e and Beltrami-Klein discs
with unit discs centred at the origin on a plane (of complex numbers).
Show that the map $z\mapsto\frac{2z}{1+|z|^2}$, up to a scale factor,
is an isometry.

\smallskip

We have forgotten to mention one more type of geodesic. It corresponds
to a subspace $W\le V$ whose hermitian form is real, nonnull, and
degenerate. In spite of the fact that the length of every segment
contained in such $\Bbb P_\Bbb KW$ vanishes, $\Bbb P_\Bbb KW$ is a {\it
bona fide\/} geodesic (see Section 4.7).

\smallskip

{\bf4.5.11.~Duality.} The hermitian form establishes a bijection
between points and geodesics in the M\"obius-Beltrami-Klein projective
plane: the point $p\in\Bbb P_\Bbb RV$ corresponds to the geodesic
$\Bbb P_\Bbb Rp^\perp$. If $p$ is negative/positive, then
$\Bbb P_\Bbb Rp^\perp$ is spherical/hyperbolic. If $p$ is isotropic,
then $\Bbb P_\Bbb Rp^\perp$ is a degenerate geodesic (with $p^\perp$ of
signature $0+$) that is tangent to the absolute and passes through $p$.

On the one hand, a hyperbolic geodesic is simply a pair of distinct
points at the absolute (its vertices). On the other hand, a hyperbolic
geodesic in the Beltrami-Klein disc is given by a positive point. This
means that the M\"obius band $\E V$ equipped with its lorentzian metric
describes the geometry of the space of geodesics in the Beltrami-Klein
disc.

\medskip

{\bf4.6.~Space of circles.} In this section, we deal with the
Riemann-Poincar\'e sphere and study the geometry of `linear' subspaces
of the form $\Bbb P_\Bbb CW$, where $W\le V$ is a $2$-dimensional
$\Bbb R$-linear subspace.

When $W$ is a $\Bbb C$-linear subspace, $\Bbb P_\Bbb CW$ is a point.
What about the other cases? We will learn that the remaining linear
subspaces $\Bbb P_\Bbb CW$ are geometrically classified by the
signature of the form $(-,-):=\Re\langle-,-\rangle$ on $W$.

\vskip30pt

\medskip

\newpage

{\bf4.7.~Complex hyperbolic zoo.}

\vskip12pt

\noindent
$
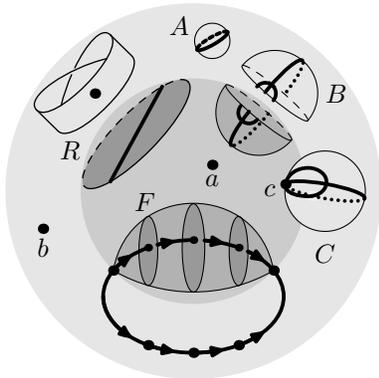
\epsfbox{Pics.10}$

\leftskip159pt

\vskip-147pt

\noindent
{\eightpoint The darker ball is the complex hyperbolic plane $\B V$ and
the lighter one is $\E V$.

\noindent
{\bf(a)}, {\bf(b)}, and {\bf(c)} are negative, positive, and isotropic
points in the extended complex hyperbolic plane. They are respectively
dual to the complex geodesics (A), (B), and (C).

\noindent
{\bf(A)} is a round sphere.\footnote{Well, with negative definite
metric.}

\noindent
{\bf(B)} is a Riemann-Poincar\'e sphere divided by the absolute into
its hyperbolic discs. Two geodesics and the absolutes are drawn.

\noindent
{\bf(C)} is a degenerate complex geodesic. Excluding the isotropic
point (c), its~geometry is affine. Two geodesics are drawn.

\noindent
{\bf(R)} is a M\"obius-Beltrami-Klein projective plane (commonly called
an $\Bbb R$-{\it plane\/}). The point and the geodesic are dual to each
other inside the plane (the extension of the geodesic to the band is
not in the picture).

\noindent
{\bf(F)} is a bisector. Its slices and real spine are drawn. Every
slice is a hyperbolic disc (complex geodesic) dual to a point in $\E$
belonging to the real spine.}

\leftskip0pt

\vskip130pt

\medskip

{\bf4.8.~Finite configurations.} In 1872, Felix Klein came up with a
brilliant idea: in geometry, one~should study the properties of a space
which are invariant under the symmetries of the space. This view became
known as the {\it Erlangen Program.} It was, and still is, very
revolutionary. Let us give some examples at the level of plane
Euclidean geometry. We are used to characterize some triangles in the
Euclidean plane as being equal\footnote{It is certain that absolute
equality does not exist in the real world. But it does not exist in the
mathematical world either! Do you mean that $1=1$ is an absolute
equality? No way! This `equality' just expresses the fact that two sets
of one element are equivalent in the sense that there exists a
bijection between them. For example, $1=1$ does not imply that one
person equals another (which seems to be very good!).}
while, in fact, they are not equal as subsets in the plane. The
triangles are {\it geometrically equal,} that is, there exists a
symmetry of the plane (a geometry-preserving bijection) that sends one
triangle onto the other. The composition of symmetries and the inverse
of a symmetry are symmetries. In other words, the symmetries constitute
a group (see Section 7 for the definition).

Roughly speaking, geometry is not made of objects, but of objects and
movements. The allowed movements vary from case to case and, generally,
we can study the geometry of any structure. {\sl This means that we
actually study the symmetry group of the structure.} A simple example:
studying the geometry of a set with no imposed structure is the study
of the permutation group of the set. A difficult example:

\medskip

\hskip130pt{\oitossi You boil it in sawdust, you salt it in glue

\hskip130pt You condense it with locusts and tape

\hskip130pt Still keeping one principal object in view ---

\hskip130pt To preserve its symmetrical shape.}

\medskip

\hskip130pt{\oitoss --- LEWIS CARROLL,} {\oitossi The Hunting of the
Snark}

\medskip

In what follows, one can find an intermediate example.

\smallskip

It is easy to figure out that the symmetries of a $\Bbb K$-linear space
$V$ equipped with a hermitian form $\langle-,-\rangle$ are all the
$\Bbb K$-linear isomorphisms $g:V\to V$ preserving $\langle-,-\rangle$.
They constitute the {\it unitary\/} group
$$\U V:=\big\{g\in\GL V\mid\langle gv,gv'\rangle=\langle
v,v'\rangle\text{ for all }v,v'\in V\big\}.$$

The Gram matrix provides the geometrical classification of generic
finite configurations in $V$ (finite configuration = finite tuple of
points) :

\smallskip

{\bf4.8.1.~Stollen Carlos' lemma.} {\sl Let\/ $w_1,w_2,\dots,w_k\in V$
and\/ $w'_1,w'_2,\dots,w'_k\in V$ be configurations such that the
subspaces\/ $W:=\Bbb Kw_1+\Bbb Kw_2+\dots+\Bbb Kw_k$ and\/
$W':=\Bbb Kw'_1+\Bbb Kw'_2+\dots+\Bbb Kw'_k$ are nondegenerate. Then
the configurations are geometrically equal, i.e., there exists\/
$g\in\U V$ such that\/ $gw_i=w'_i$ for all\/ $i$, iff their Gram
matrices\/ $G(w_1,w_2,\dots,w_k)$ and\/ $G(w'_1,w'_2,\dots,w'_k)$
are equal.}

\smallskip

{\bf Proof.}~If such a $g$ exists, then
$\langle w'_i,w'_j\rangle=\langle gw_i,gw_j\rangle=\langle
w_i,w_j\rangle$
for all $i,j$ since $g\in\U V$. In other words,
$G(w'_1,w'_2,\dots,w'_k)=G(w_1,w_2,\dots,w_k)$.

Conversely, suppose that
$G(w_1,w_2,\dots,w_k)=G(w'_1,w'_2,\dots,w'_k)$. We define the linear
map $h:\Bbb K^k\to\nomathbreak W$,
$h:(c_1,c_2,\dots,c_k)\mapsto\sum\limits_{i=1}^kc_iw_i$. Obviously, $h$
is surjective. In a similar way, we define the surjective linear map
$h':\Bbb K^k\to W'$. Let us prove that $\ker h=\ker h'$. By symmetry,
it suffices to show that $\ker h\subset\ker h'$. If
$(c_1,c_2,\dots,c_k)\in\ker h$, that is, if
$\sum\limits_{i=1}^kc_iw_i=0$, then
$$0=\Big\langle\sum\limits_{i=1}^kc_iw_i,w_j\Big\rangle=
\sum\limits_{i=1}^kc_i\langle w_i,w_j\rangle=
\sum\limits_{i=1}^kc_i\langle w'_i,w'_j\rangle=
\Big\langle\sum\limits_{i=1}^kc_iw'_i,w'_j\Big\rangle$$
for all $j$. Being $W'$ nondegenerate, we have
$\sum\limits_{i=1}^kc_iw'_i=0$, that is,
$(c_1,c_2,\dots,c_k)\in\ker h'$.

We obtained a linear isomorphism $l:W\to W'$ such that $lw_i=w'_i$ for
all $i$. It follows from $G(w_1,w_2,\dots,w_k)=G(w'_1,w'_2,\dots,w'_k)$
that $l$ preserves the form, that is,
$\langle lx,ly\rangle=\langle x,y\rangle$ for all $x,y\in W$. In
particular, $W$ and $W'$ are of the same signature. By Exercise 6.6, we
have orthogonal decompositions $V=W\oplus W^\perp$ and
$V=W'\oplus{W'}^\perp$. Hence, $W^\perp$ and ${W'}^\perp$ are of the
same signature. Therefore, there exists a linear isomorphism
$l':W^\perp\to{W'}^\perp$ that preserves the form. It remains to define
$g:V\to V$ by the rule $g:w+u\mapsto lw+l'u$, where $w\in W$ and
$u\in W^\perp$
$_\blacksquare$

\smallskip

{\bf4.8.2.*~Exercise.}~Find necessary and sufficient conditions for the
geometric equality of two finite configurations without the assumption
that $W$ and $W'$ are nondegenerate.

\vskip40pt

\medskip

\newpage

{\bf4.9.~There is no sin south of the equator.\footnote{A quote from
the famous brazilian musician Chico Buarque.}} The word trigonometry
stands in the Greek for `measuring triangles.' The typical approach to
studying triangles in non-Euclidean plane geometry is to write down
several identities that relate, via trigonometric and hyperbolic
trigonometric functions like $\sin$, $\cos$, $\sinh$, $\cosh$, etc.,
the angles and the lengths of the sides of a triangle. Since high
school, we are used to `solve' triangles via trigonometry $\dots$ let
us see how the study of finite configurations in classic geometries
developed in the previous section may help in understanding where
trigonometric relations come from.

We begin with spherical plane geometry. As in the first of Examples
4.4, let $V$ be a $2$-dimensional complex linear space with a hermitian
form $\langle-,-\rangle$ of signature $++$. The Riemann sphere
$\Bbb P_\Bbb CV$ endowed with the metric (4.3.1) is the round
sphere of radius $\frac12$. Let $p_1,p_2,p_3\in\Bbb P_\Bbb CV$ be
distinct points such that $\langle p_i,p_j\rangle\ne0$ for all $i,j$.
They determine the oriented triangle $\Delta(p_1,p_2,p_3)$ whose side
$p_ip_{i+1}$ is the shortest segment of geodesic joining $p_i$ and
$p_{i+1}$ (the indices are modulo $3$). In particular,
$l_i:=\ell(p_ip_{i+1})<\frac\pi2$. We know from Exercise 4.8.? that
there exist representatives $p_1,p_2,p_3\in V$ with the Gram matrix
$\left[\smallmatrix1&r_1&r_3\overline\varepsilon\\r_1&1&r_2\\r_3
\varepsilon&r_2&1\endsmallmatrix\right]$,
where $0<r_i<1$ and $\varepsilon\in\Bbb C$ with $|\varepsilon|=1$. The
geometrical meaning of every number in this matrix is known:
$r_i=\sqrt{\ta(p_i,p_{i+1})}$ and
$\arg\varepsilon=2\area\Delta(p_1,p_2,p_3)$. So, the $r_i$'s speak of
the lengths of the sides of $\Delta(p_1,p_2,p_3)$ while $\varepsilon$
provides the oriented area of the triangle. Being $p_1,p_2,p_3$
linearly dependent, the determinant of the Gram matrix vanishes:
$$1+2r_1r_2r_3\Re\varepsilon-r_1^2-r_2^2-r_3^2=0.
\leqno{\bold{(4.9.1)}}$$
This equation is the only relation between the geometric invariants
$r_1,r_2,r_3,\varepsilon$ (not counting inequalities). This is the {\it
fundamental trigonometric\/} identity, and any other one is derivable
from it!

For instance, the {\it first law\/} of {\it cosines\/} in spherical
trigonometry states that
$$\cos(2l_3)=\cos(2l_1)\cos(2l_2)+\cos\alpha\sin(2l_1)\sin(2l_2)$$
under the condition $0<\alpha<\pi$ for the interior angle $\alpha$ at
$p_2$. In order to deduce this law from (4.9.1), we remind the relation
between length and tance in the spherical geometry: $l_i=\arccos r_i$
(see~Subsection~4.5.4). It follows that $\cos(2l_i)=2r_i^2-1$ and
$\sin(2l_i)=2r_i\sqrt{1-r_i^2}$. So, the first law of cosines is
equivalent to
$$\cos\alpha=\frac{r_1^2+r_2^2+r_3^2-2r_1^2r_2^2-1}{2r_1r_2
\sqrt{1-r_1^2}\cdot\sqrt{1-r_2^2}}.\leqno{\bold{(4.9.2)}}$$
By Exercise 4.5.3, the tangent vectors
$$t_1:=\langle-,p_2\rangle\frac{\pi[p_2]p_1}{\langle
p_1,p_2\rangle},\qquad
t_2:=\langle-,p_2\rangle\frac{\pi[p_2]p_3}{\langle p_3,p_2\rangle}$$
are respectively tangent to $p_2p_1$ and $p_2p_3$ at $p_2$.
Therefore,
$$\cos\alpha=\frac{\Re\langle t_1,t_2\rangle}{\sqrt{\langle
t_1,t_1\rangle}\cdot\sqrt{\langle t_2,t_2\rangle}}=
\frac{r_3\Re\varepsilon-r_1r_2}{\sqrt{1-r_1^2}\cdot\sqrt{1-r_2^2}}.$$
Using the fundamental trigonometric identity (4.9.1), it is easy to see
that the above expression is exactly (4.9.2).

\smallskip

{\bf4.9.3.~Exercise.} Derive the {\it law\/} of {\it sines\/}
$$\frac{\sin(2l_1)}{\sin\alpha_3}=\frac{\sin(2l_2)}{\sin\alpha_1}=
\frac{\sin(2l_3)}{\sin\alpha_2}$$
in spherical plane geometry assuming that the length $l_i$ of the side
$p_ip_{i+1}$ and the interior angle $\alpha_i$ at the vertex $p_i$ of
the triangle $\Delta(p_1,p_2,p_3)$ satisfy the inequalities
$0<l_i<\frac\pi2$ and $0<\alpha_i<\pi$, $i=1,2,3$.

\smallskip

{\bf4.9.4.~Exercise.} Let $\Delta(p_1,p_2,p_3)$ be a triangle in the
Riemann-Poincar\'e sphere with distinct nonisotropic vertices of the
same signature. Write down the fundamental trigonometric identity for
the triangle and derive the {\it first\/} and {\it second laws\/} of
{\it cosines\/} as well as the {\it law\/} of {\it sines\/} in
hyperbolic geometry:
$$\cosh(2l_3)=\cosh(2l_1)\cosh(2l_2)-
\cos\alpha_2\sinh(2l_1)\sinh(2l_2),$$
$$\cos\alpha_2+\cos\alpha_2\cos\alpha_3=
\cosh(2l_3)\sin\alpha_2\sin\alpha_3,$$
$$\frac{\sinh(2l_1)}{\sin\alpha_3}=\frac{\sinh(2l_2)}{\sin\alpha_1}=
\frac{\sinh(2l_3)}{\sin\alpha_2},$$
where $l_i$ stands for the length of the side $p_ip_{i+1}$ and
$\alpha_i$, for the interior angle at $p_i$ (the indices are modulo
$3$). Study the trigonometry of triangles with the other signatures of
vertices (including isotropic ones).

\medskip

{\bf4.10A.~Geometry on the absolute.}

\vskip40pt

\medskip

\newpage

{\bf4.11.~A bit of history.} In 1820, the eighteen years old hungarian
mathematician J\'anos Bolyai began to write a treatise on non-Euclidean
geometry. His father, Farkas Bolyai, had himself struggled in vain with
the parallel postulate for many years. Farkas Bolyai did not measure
efforts in trying to dissuade his son from following what he thought
was a hopeless path:

\smallskip

`{\sl You must not attempt this approach to parallels. I know this way
to its very end. I have traversed this bottomless night, which
extinguished all light and joy of my life. I entreat you, leave the
science of parallels alone $\dots$ I thought I would sacrifice myself
for the sake of truth. I was ready to become a martyr who would remove
the flaw from geometry and return it purified to mankind. I
accomplished monstrous, enormous labors\/{\rm;} my creations are far
better than those of others and yet I have not achieved complete
satisfaction $\dots$ I turned back when I saw that no man can reach the
bottom of the night. I~turned back unconsoled, pitying myself and all
mankind.

I admit that I expect little from the deviation of your lines. It seems
to me that I have been in these regions\/{\rm;} that I have traveled
past all reefs of this infernal Dead Sea and have always come back with
broken mast and torn sail. The ruin of my disposition and my fall date
back to this time. I thoughtlessly risked my life and happiness ---}
{\it aut Caesar aut nihil.}'

\smallskip

\noindent
Yet, J\'anos had enough courage to pursuit his ideas. And where many
failed, the young genius succeeded~$\dots$ He wrote to his father:

\smallskip

`{\sl It is now my definite plan to publish a work on parallels as soon
as I can complete and arrange the material and an opportunity presents
itself $\dots$ I have discovered such wonderful things that I was
amazed, and it would be an everlasting piece of bad fortune if they
were lost. When you, my dear Father, see them, you will
understand\/{\rm;} at present I can say nothing except this\/{\rm:}
that out of nothing I have created a strange new universe. All that I
have sent you previously is like a house of cards in comparison with a
tower. I am no less convinced that these discoveries will bring me
honor than I would be if they were completed.}'

\smallskip

Naturally, J\'anos desired to present his discoveries to the foremost
of the mathematicians, the {\it princeps mathematicorum,} Carl
Friedrich Gauss. It turns out that Farkas Bolyai was old friends with
Gauss, and~the opportunity J\'anos was so eagerly looking for stood
right in front of him: his father would write a letter to Gauss and
communicate his son's great accomplishments. It could not get any
better.

\smallskip

Finally, an answer from Gauss to Farkas arrived:

\smallskip

`{\sl If I begin with the statement that I dare not praise such a work,
you will of course be startled for a moment\/}:'

\smallskip

Why not praise my work? --- thought J\'anos. Is it possible that
everything is wrong? Have I, like many, fell in some of the elusive
traps surrounding the parallels? No, it must not be!

\smallskip

\noindent
`{\sl but I cannot do otherwise\/};' --- proceeded Gauss --- `{\sl to
praise it would amount to praising myself\/{\rm;} for~the entire
content of the work, the path which your son has taken, the results to
which he is led, coincide almost exactly with my own meditations which
have occupied my mind for from thirty to thirty-five years.}'

\smallskip

That definitely was not fair! --- thought J\'anos --- Could not Gauss
acknowledge honestly, definitely, and frankly my work? Verily, it is
not this attitude we call life, work, and merit. J\'anos was so
profoundly disappointed that he could never fully recover from this
episode.

\smallskip

It was a rainy evening,  October 17 1841, when J\'anos received from
his father a brochure entitled, to~his surprise, `Geometrische
Untersuchengen zur Theorie der Parallellinien' (Geometrical
investigations on the theory of parallel lines). J\'anos was a polyglot
and spoke perfectly nine foreign languages. Reading German was no
challenge to him. The author of the brochure? Some russian professor
Nikolai Ivanovich Lobachevsky. The more J\'anos read the brochure, the
more puzzled he got. All his cherish discoveries, the great discoveries
no one would ever acknowledge him for, they were all there $\dots$ he
flipped the pages with more and more anguish $\dots$ no doubts, the
work in his hands was a masterpiece. J\'anos closed the brochure and
left it on the table. He took a few steps back and just glanced at the
brochure for a while. A russian professor Nikolai Ivanovich Lobachevsky
$\dots$ that writes a beautiful text in German about non-Euclidean
geometry $\dots$ J\'anos eyes became injected with rage and he punched
the table furiously. This is the last straw! --- he cried. J\'anos
utmost suspicion, naturally, was that no professor Lobachevsky ever
existed, and that the brochure was nothing but a work of Gauss.

\smallskip

{\bf4.11.1.~References.} For the correspondence between Farkas Bolyai
and J\'anos Bolyai see [Mes]. For~the story about the brochure see
[Kag, p.~391, l.~13--15]. For Gauss' correspondence see [Sch].

\smallskip

\noindent
[Kag] Kagan, V.~F., {\it Lobachevsky,} edition of the Academy of
Sciences of the USSR, Moscow-Leningrad, 1948 (Russian)

\smallskip

\noindent
[Mes] Meschkowski, H., {\it Evolution of mathematical thought,}
Holden-Day, San Francisco, 1965

\smallskip

\noindent
[Sch] Schmidt, F., and St\"ackel, P., {\it Briefwechsel zwischen
C.~F.~Gauss and W.~Bolyai,} Johnson Reprint Corp. New York, 1972
(German)

\rightheadtext{Riemann surfaces}

\bigskip

\centerline{\bf5.~Riemann surfaces}

\medskip

{\bf5.1.~Regular covering and fundamental group.}

\medskip

{\bf5.2.~Discrete groups and Poincar\'e polygonal theorem.}

\medskip

{\bf5.3.~Teichm\"uller space.}

\rightheadtext{Appendix: Largo al factotum della citta}

\newpage

\bigskip

\centerline{\bf6.~Appendix: Largo al factotum della citta}

\medskip

If you leave the Universidad de Sevilla and walk down the Calle Palos
de la Frontera street (heading the Plaza de Espa\~na), you might
unexpectedly hear the melody

\medskip

\hskip180pt{\sl Rasori e pettini}

\hskip180pt{\sl lancette e forbici}

\hskip180pt{\sl al mio comando}

\hskip180pt{\sl tutto qui sta.\footnote{Raisors and combs, blades and
scissors at my disposal here they are.}}

\medskip

\noindent
coming out of a barber shop. It sounds so familiar that you decide to
enter the shop. The barber introduces himself:

\medskip

--- Ciao, mi chiamo Figaro, il barbiere-factotum.\footnote{Hello, I am
Figaro, a factotum barber.}

--- Hi, I am a student of mathematics here at the University.

--- Hum, a mathematician $\dots$ The mathematicians use to look for me
only for two reasons $\dots$ --- Figaro seems annoyed.

--- $\dots$ they do not know how to solve the Barber
Paradox\footnote{Also known as Russell's Paradox (Bertrand Russell,
British philosopher and mathematician) : Who shaves the barber that
shaves only men that do not shave themselves?}
$\dots$

--- $\dots$ or they cannot solve their problems because they do not
know the linear tools, like Linear Algebra! To say nothing of Hermitian
Tools! --- Figaro is now furious.

\medskip

\noindent
You may become confused. It is a comprehensible thing that
mathematicians could seek the barber to get convinced of his existence.
But $\dots$

\medskip

--- Why on earth would an ignorant in Linear Algebra look for you?

--- Not knowing Linear Algebra is a barbarity. And I am a barber, what
do you expect? Sit down and let me introduce to you the linear and
hermitian tools:

\medskip

\hskip180pt{\sl Rasori e pettini}

\hskip180pt{\sl lancette e forbici}

\hskip180pt{\sl al mio comando}

\hskip180pt{\sl tutto qui sta.}

\medskip

We deal with finite-dimensional linear spaces over $\Bbb R$ or
$\Bbb C$. To cover both cases, denote the scalars by $\Bbb K$. The
symbol $\overline k$ stands for the {\it conjugate\/} to the (complex)
number $k\in\Bbb K$.

\medskip

{\bf6.1.~Definition.} Let $V$ be a $\Bbb K$-linear space. A {\it
hermitian form\/} is a map $\langle-,-\rangle:V\times V\to\Bbb K$,
$(x,y)\mapsto\langle x,y\rangle$ linear in $x$ and such that
$\langle x,y\rangle=\overline{\langle y,x\rangle}$ for all $x,y\in V$.
In other words, the form is $1.5$-{\it linear\/} since
$\langle kx,y\rangle=k\langle x,y\rangle$ and
$\langle x,ky\rangle=\overline k\langle x,y\rangle$ for all
$k\in\Bbb K$. If $W\le V$ is a subspace, then we can restrict the form
$\langle-,-\rangle$ to $W$, getting a linear space $W$ equipped with
the {\it induced\/} hermitian form.

\medskip

{\bf6.2.~Definition.} Let $V$ be a linear space equipped with a
hermitian form and let $W\le V$ be a subspace. We define
$W^\perp:=\big\{v\in V\mid\langle v,W\rangle=0\big\}$, the {\it
orthogonal\/} to $W$. We call $V^\perp$ the {\it kernel\/} of the form
on $V$. If the kernel vanishes, we say that the form is {\it
nondegenerate.} If the induced form on a subspace $W\le V$ is
nondegenerate, $W$ is said to be {\it nondegenerate.} For $U,W\le V$,
the {\it orthogonal\/} of $W$ {\it relatively\/} to $U$ is given by
$W^{\perp U}:=W^\perp\cap U$.

\medskip

{\bf6.3.~Exercise.} Show that $W^\perp\le V$ and
$W\subset{W^\perp}^\perp$ for all $W\le V$. Prove also that
$(W_1+W_2)^\perp=W_1^\perp\cap W_2^\perp$ for all $W_1,W_2\le V$. Is
the identity $(W_1\cap W_2)^\perp=W_1^\perp+W_2^\perp$ true?

\medskip

{\bf6.4.~Exercise.} Define the induced form on $V/V^\perp$ and verify
that this definition is correct. Show that $V/V^\perp$ is
nondegenerate. Decomposing $V=V^\perp\oplus W$, prove that the spaces
$V/V^\perp$ and $W$ equipped with the induced forms are naturally
isomorphic.

\medskip

{\bf6.5.~Exercise.} For $W\le V$, show that
$\dim W+\dim W^\perp\ge\dim V$.

\medskip

{\bf6.6.~Exercise.} Show that $V=W\oplus W^\perp$ for every
nondegenerate subspace $W\le V$.

\medskip

{\bf6.7.~Exercise.} Suppose that both $W$ and $V$ are nondegenerate,
where $W\le V$. Prove that ${W^\perp}^\perp=\nomathbreak W$.
\medskip

{\bf6.8.~Exercise.} Suppose that both $W$ and $V$ are nondegenerate,
where $W\le V$. Show that $W^\perp$ is nondegenerate.

\medskip

{\bf6.9.~Exercise.} Show that there exists a {\it nonisotropic\/}
$v\in V$, i.e., $\langle v,v\rangle\ne0$, if
$\langle-,-\rangle\not\equiv0$.

\medskip

{\bf6.10.~Exercise.} Suppose that both $W$ and $V$ are nondegenerate,
where $W\lvertneqq V$. Show that there exists a nondegenerate subspace
$W'\le V$ such that $W\le W'$ and $\dim W'=\dim W+1$.

\medskip

{\bf6.11.~Definition.} A {\it flag\/} of subspaces is a chain of
subspaces $V_0\le V_1\le\dots\le V_n$ such that $V_n=V$ and
$\dim V_i=i$ for all $i$. If $V$ is equipped with a hermitian form, a
flag is {\it nondegenerate\/} when all $V_i$'s are nondegenerate.

\medskip

{\bf6.12.~Exercise.} Show that every nondegenerate linear space admits
a nondegenerate flag of subspaces.

\medskip

{\bf6.13.~Definition.} A linear basis $\beta:b_1,b_2,\dots,b_n$ is {\it
orthonormal\/} if $\langle b_i,b_i\rangle\in\{-1,0,1\}$ and
$\langle b_i,b_j\rangle=0$ for all $i$ and $j$ such that $i\ne j$.
Denote by $\beta_-,\beta_0,\beta_+$ the amount of elements in the basis
$\varepsilon$ such that $\langle b_i,b_i\rangle=-1$,
$\langle b_i,b_i\rangle=0$, $\langle b_i,b_i\rangle=1$, respectively.
The triple $(\beta_-,\beta_0,\beta_+)$ is the {\it signature\/} of the
basis.

\medskip

{\bf6.14.~Exercise.} Let $\beta:b_1,b_2,\dots,b_n$ be an orthonormal
basis in $V$. Show that $\beta_0$ is the dimension of the kernel of the
form on $V$, $\beta_0=\dim V^\perp$.

\medskip

{\bf6.15.~Gram-Schmidt orthogonalization.} {\sl Let\/
$V_0\le V_1\le\dots\le V_n$ be a nondegenerate flag of subspaces in\/
$V$. Then there exists an orthonormal basis\/ $b_1,b_2,\dots,b_n$ in\/
$V$ such that\/ $b_1,b_2,\dots,b_k$ is a basis in\/ $V_k$ for all\/
$k$.}

\medskip

{\bf Proof.} Induction on $n$. For $n=1$, we simply take some
$0\ne c_1\in V_1$ and normalize it:
$b_1:=\displaystyle\frac{c_1}{\sqrt{\big|\langle
c_1,c_1\rangle\big|}}$.
(Being $V_1$ nondegenerate, $\langle c_1,c_1\rangle\ne0$.) Suppose
that, for some $k<n$, we have already found an orthonormal basis
$b_1,b_2,\dots,b_k$ in $V_k$ such that $b_1,b_2,\dots,b_i$ is a basis
in $V_i$ for all $i\le k$. We~choose $c_{k+1}\in V_{k+1}\setminus V_k$
and put
$c'_{k+1}=c_{k+1}-\displaystyle\sum\limits_{i=1}^k\frac{\langle
c_{k+1},b_i\rangle}{\langle b_i,b_i\rangle}b_i$.
Taking into account that the $b_i$'s are orthogonal, a straightforward
calculus shows that $\langle c'_{k+1},b_i\rangle=0$ for all $i\le k$.
If $c'_{k+1}$ could be isotropic, then it would belong to the kernel of
the form on $V_{k+1}$. Therefore, $c'_{k+1}$ is nonisotropic and we can
normalize $c'_{k+1}$, getting the desired $b_{k+1}$
$_\blacksquare$

\medskip

{\bf6.16.~Corollary.} {\sl Every linear space with a hermitian form
admits an orthonormal basis.}

\medskip

{\bf Proof.}~By Exercise 6.4, we can assume that the space $V$ is
nondegenerate. Using Exercise 6.10, we~can build a nondegenerate flag
of subspaces in $V$. Now, the result follows from 6.15
$_\blacksquare$

\medskip

{\bf6.17.~Definition.} Let $v_1,v_2,\dots,v_k\in V$. The matrix
$G:=G(v_1,v_2,\dots,v_k):=[g_{ij}]$, where
$g_{ij}:=\langle v_i,v_j\rangle$, is called the {\it Gram\/} matrix of
$v_1,v_2,\dots,v_k$.

\medskip

Obviously, $\overline G^t=G$, where $M^t$ denotes the transpose matrix
of $M$ and $\overline M$ denotes the matrix $M$ with conjugate entries.
In other words, $G$ is {\it hermitian\/} ({\it symmetric\/}).

The Gram matrix $G^{\beta\beta}:=G(b_1,b_2,\dots,b_n)$ of some basis
$\beta:b_1,b_2,\dots,b_n$ in $V$ determines the hermitian form on $V$
since
$\langle v,v'\rangle=[v]_\beta^tG^{\beta\beta}\overline{[v']}_\beta$
for all $v,v'\in V$, where $[v]_\beta$ denotes the column matrix whose
entries are the coefficients $c_i$ appearing in the linear combination
$v=\sum\limits_{i=1}^nc_ib_i$. Indeed, if~$v=\sum\limits_{i=1}^nc_ib_i$
and $v'=\sum\limits_{i=1}^nc'_ib_i$, then
$\langle v,v'\rangle=\sum\limits_{i,j=1}^nc_i\langle
b_i,b_j\rangle\overline{c'}_j=
\sum\limits_{i,j=1}^nc_ig_{ij}\overline{c'}_j$.
A basis is orthonormal iff its Gram matrix is diagonal with diagonal
entries $-1,0,1$. We emphasize that every hermitian matrix is the Gram
matrix of a basis in a certain linear space with an appropriate
hermitian form.

Let $\alpha:a_1,a_2,\dots,a_n$ be another basis in $V$ and let
$M^\beta_{\alpha}=[m_{ij}]$ be the matrix representing a change of
basis from $\alpha$ to $\beta$, that is,
$b_j=\sum\limits_{i=1}^nm_{ij}a_i$ for all $j$. Then
$$g_{kl}=\langle b_k,b_l\rangle=\Big\langle\sum
\limits_{i=1}^nm_{ik}a_i,\sum\limits_{i=1}^nm_{jl}a_j\Big\rangle=
\sum\limits_{i,j=1}^nm_{ik}\langle a_i,a_j\rangle\overline
m_{jl}=\sum\limits_{i,j=1}^nm_{ik}f_{ij}\overline m_{jl},$$
where $G^{\alpha\alpha}=[f_{ij}]$. We obtained the relation
$G^{\beta\beta}=
(M^\beta_{\alpha})^tG^{\alpha\alpha}\overline{M^\beta_{\alpha}}$.
In particular, it follows that the sign of $\det G^{\beta\beta}$ does
not depend on the choice of the basis because
$$\det G^{\beta\beta}=\det(M^\beta_{\alpha})^t\det
G^{\alpha\alpha}\det\overline{M^\beta_{\alpha}}=\det
M^\beta_{\alpha}\det G^{\alpha\alpha}\overline{\det
M^\beta_{\alpha}}=|\det M^\beta_{\alpha}|^2\det G^{\alpha\alpha}.$$

\medskip

{\bf6.18.~Lemma.} {\sl Let\/ $G^{\beta\beta}$ be the Gram matrix of
a basis in a linear space\/ $V$. Then\/ $V$ is degenerate iff\/
$\det G^{\beta\beta}=0$
$_\blacksquare$}

\medskip

{\bf6.19.~Example.} Let $V\ni e,f$ be such that
$\langle e,e\rangle>0>\langle f,f\rangle$. We put $W:=\Bbb Ke+\Bbb Kf$.
Then $\dim W=2$ and every orthonormal basis in $W$ has signature
$(1,0,1)$. Moreover, $W$ contains (non-null) nonisotropic elements.

Indeed, we can take $W=V$. If $0\ne n\in V^\perp$, then
$V=\Bbb Kb+\Bbb Kn$ for some $b\in V$. Assuming
$\langle b,b\rangle\ge0$ we obtain $\langle v,v\rangle\ge0$ for all
$v\in V$ and assuming $\langle b,b\rangle\le0$ we obtain
$\langle v,v\rangle\le0$ for all $v\in V$. Both cases are impossible
since $V$ contains one positive element and one negative element. For a
similar reason, $\dim V=2$. Taking an orthonormal basis $\beta$ in $V$,
it is easy to see that the signature of such basis is distinct from
$(2,0,0)$ (since $V$ contains a positive element) and from $(0,0,2)$
(since $V$ contains a negative element). By Exercise 6.14, $\beta_0=0$.
Hence, the signature is $(1,0,1)$. Obviously, the sum of the elements
of the orthonormal basis is isotropic.

\medskip

{\bf6.20.~Sylvester's law of inertia.} {\sl The signature does not
depend on the choice of an orthonormal basis.}

\medskip

{\bf Proof.} Induction on $\dim V$. By Exercises 6.4 and 6.14, we can
assume that $V$ is nondegenerate. Let $\beta:b_1,b_2,\dots,b_n$ and
$\beta':b'_1,b'_2,\dots,b'_n$ be orthonormal bases. So,
$\beta_0=\beta'_0=0$ by Exercise 6.14. If~$\beta_-=0$, then
$\langle v,v\rangle\ge0$ for all $v\in V$, implying $\beta'_-=0$. In
the same way, $\beta_+=0$ implies $\beta'_+=0$. Therefore, we can
assume that $\langle b_n,b_n\rangle=1$ and
$\langle b'_n,b'_n\rangle=-1$. We put
$$W:=\Bbb Kb_n+\Bbb Kb'_n,\qquad U:=(\Bbb Kb_n)^\perp,\qquad
U'=(\Bbb Kb'_n)^\perp.$$
It is easy to see that $U=\Bbb Kb_1+\Bbb Kb_2+\dots+\Bbb Kb_{n-1}$ and
$U'=\Bbb Kb'_1+\Bbb Kb'_2+\dots+\Bbb Kb'_{n-1}$. Therefore,
the~signatures of the indicated bases in $U$ and $U'$ are respectively
$(\beta_-,0,\beta_+-1)$ and $(\beta'_--1,0,\beta'_+)$. By Exercise 6.3,
$W^\perp=U\cap U'$. By Example 6.19 and Exercise 6.14, $W$ is
nondegenerate. So, $U\cap U'$ is nondegenerate by Exercise 6.8.
Applying Exercise 2.6 to the spaces $U$ and $U'$ and to the subspace
$U\cap U'$, we obtain the orthogonal decompositions
$U=(U\cap U')\oplus(U\cap U')^{\perp U}$ and
$U'=(U\cap U')\oplus(U\cap U')^{\perp U'}$. Using Corollary 6.16, we
choose an orthonormal basis $\alpha$ in $U\cap U'$. Let $\gamma$ and
$\gamma'$ be some orthonormal bases respectively in
$(U\cap U')^{\perp U}$ and $(U\cap U')^{\perp U'}$. Therefore,
$\alpha\sqcup\gamma$ and $\alpha\sqcup\gamma'$ are orthonormal bases
respectively in $U$ and $U'$. Calculating the signatures, we obtain
$$(\beta_-,0,\beta_+-1)=\big((\alpha\sqcup\gamma)_-,
(\alpha\sqcup\gamma)_0,(\alpha\sqcup\gamma)_+\big)=
(\alpha_-,\alpha_0,\alpha_+)+(\gamma_-,\gamma_0,\gamma_+),$$
$$(\beta'_--1,0,\beta'_+)=\big((\alpha\sqcup\gamma')_-,
(\alpha\sqcup\gamma')_0,(\alpha\sqcup\gamma')_+\big)=
(\alpha_-,\alpha_0,\alpha_+)+(\gamma'_-,\gamma'_0,\gamma'_+)$$
since the signatures do not depend on the choices of orthogonal bases
in $U$ and $U'$ by the induction hypothesis. It remains to show that
$(U\cap U')^{\perp U}=(\Bbb Kb_n)^{\perp W}$ and that
$(U\cap U')^{\perp U'}=(\Bbb Kb'_n)^{\perp W}$ since this implies
$(\gamma_-,\gamma_0,\gamma_+)=(1,0,0)$ and
$(\gamma'_-,\gamma'_0,\gamma'_+)=(0,0,1)$ by Example 2.19.

Being $W$ and $V$ nondegenerate,
$(U\cap U')^{\perp U}=(U\cap U')^\perp\cap U={W^\perp}^\perp\cap
U=W\cap(\Bbb Kb_n)^\perp=(\Bbb Kb_n)^{\perp W}$
by Exercise 6.7. For the same reason,
$(U\cap U')^{\perp U'}=(\Bbb Kb'_n)^{\perp W}$
$_\blacksquare$

\medskip

We can now speak of the {\it signature\/} of a space. How do we measure
it? By Exercise 6.14, $\beta_0=\dim V^\perp$. Using Exercise 6.4, the
problem can be reduced to the case of a nondegenerate $V$. Let
$\gamma:c_1,c_2,\dots,c_n$ be a basis in $V$ with a known Gram matrix
$G^{\gamma\gamma}$. We want to find out the signature of $V$ in terms
of $G^{\gamma\gamma}$. Defining
$V_k:=\Bbb Kc_1+\Bbb Kc_2+\dots+\Bbb Kc_k$ for every $0\le k\le n$, we
obtain a flag of subspaces. Obviously, the Gram matrix of the basis
$c_1,c_2,\dots,c_k$ in $V_k$ is the $(k\times k)$-submatrix
$G_k^{\gamma\gamma}$ (called a {\it principal\/} submatrix) formed by
the first $k$ lines and by the first $k$ columns of
$G^{\gamma\gamma}=G_n^{\gamma\gamma}$. We assume that the flag is
nondegenerate. By Lemma 6.18, this is equivalent to
$\det G_k^{\gamma\gamma}\ne0$ for all $1\le k\le n$. We
apply\footnote{As it usually happens, the proof is more important than
the fact itself.}
Orthogonalization 6.15 to the flag and observe that the signs of the
determinants $\det G_k^{\gamma\gamma}$ related to the bases
$b_1,b_2,\dots,b_k,c_{k+1},c_{k+2},\dots,c_n$ do not change when we
increase $k$ because the first $l$ elements in
$b_1,b_2,\dots,b_k,c_{k+1},c_{k+2},\dots,c_n$ constitute a basis in
$V_l$ for all $l$. When we arrive at an orthonormal basis, the
signature can be measured as follows:

\medskip

{\bf6.21.~Sylvester criterion.}~{\sl If\/ $\det G_k^{\gamma\gamma}\ne0$
for every\/ $1\le k\le n$, then the signature of the space equals\/
$(n_-,0,n_+)$, where\/ $n_-$ is the amount of negative numbers in the
sequence
$$\det G_1^{\gamma\gamma},\qquad\frac{\det G_2^{\gamma\gamma}}{\det
G_1^{\gamma\gamma}},\qquad\frac{\det G_3^{\gamma\gamma}}{\det
G_2^{\gamma\gamma}},\qquad\dots,\qquad\frac{\det G_n^{\gamma\gamma}}{\det
G_{n-1}^{\gamma\gamma}}$$
and\/ $n_+$ is the amount of positive numbers in the same sequence
$_\blacksquare$}

\medskip

{\bf6.22.*~Exercise.}~Find a criterion without the assumption that
$\det G_k^{\gamma\gamma}\ne0$ for every $k$.

\medskip

Exercises 6.23--26 concern the study of the possible signatures of a
subspace when the signature of the space is given. Note that two spaces
of the same signature admit an isomorphism between them that preserves
the form.

\rightheadtext{Appendix: Basic algebra and topology}

\medskip

{\bf6.23.~Exercise.} Let $V$ be a space of signature $(n_-,n_0,n_+)$.
Show that $V$ contains a subspace $W$ of signature $(m_-,m_0,m_+)$ iff
the space $V/V^\perp$ (of signature $(n_-,0,n_+)$) possesses a subspace
of signature $(m_-,m_0-m,m_+)$ for some $m$ such that $0\le m\le n_0$.

\medskip

{\bf6.24.~Exercise.} Let $V$ be a space of signature $(n_-,0,n_+)$.
Show that $\min(n_-,n_+)$ is the highest possible dimension of a
subspace $W$ with the null induced form.

\medskip

{\bf6.25.~Exercise.} Let $V$ be a space of signature $(n_-,0,n_+)$.
Show that $V$ contains a subspace $W$ of signature $(m_-,m_0,m_+)$ iff
$$m_-\le n_-,\qquad m_+\le n_+,\qquad m_0\le n_--m_-,\qquad m_0\le
n_+-m_+.$$

{\bf6.26.~Exercise.} Let $V$ be a space of signature $(n_-,n_0,n_+)$.
Show that $V$ contains a subspace of signature $(m_-,m_0,m_+)$ iff
$$m_-\le n_-,\qquad m_+\le n_+,\qquad m_-+m_0\le n_-+n_0,\qquad
m_0+m_+\le n_0+n_+.$$

\bigskip

\centerline{\bf7.~Appendix: Basic algebra and topology}

\medskip

\bigskip

\centerline{\bf8.~Appendix: Classification of compact surfaces}

\medskip

\bigskip

\centerline{\bf9A.~Appendix: Riemannian geometry}

\medskip

\bigskip

\centerline{\bf10A.~Appendix: Hyperelliptic surfaces and Goldman's
theorem}

\medskip

\bigskip

\newpage

\rightheadtext{Hints}

\centerline{\bf Hints}

\medskip

{\bf1.2.}~The question makes no sense.

\smallskip

{\bf2.3.}~Draw two distinct lines $L_1,L_2$ passing through $p$ that
are not parallel to $R_1,R_2$ and denote the intersections
$\{q_{ij}\}=R_i\cap L_j$. Joining $q_{11},q_{22}$ and $q_{12},q_{21}$,
we respectively obtain the lines $D_1$ and $D_2$. They intersect in
$\{d\}=D_1\cap D_2$. Denoting $\{q_i\}=R_i\cap L$, where $L$ is the
line joining $p$ and $d$, we can construct the lines $S_1$ and $S_2$
that join respectively $q_1,q_{22}$ and $q_{11},q_2$. We claim that the
intersection $\{q\}=S_1\cap S_2$ lives in the desired line $R$. To
prove this fact, choose the line joining $p$ and $b$ as the infinity,
where $\{b\}=R_1\cap R_2$.

\smallskip

{\bf2.5.}~Return to this exercise after studying the Beltrami-Klein
plane (see 4.5.11).

\smallskip

{\bf2.10.}~By induction on dimension, it suffices to deal with
subspheres of codimension $1$. Such a subsphere can be described as
$S:=\{q\in\Bbb S^n\mid fq=\varepsilon\}$, where
$0\ne f\in V^*:=\Lin_\Bbb R(V,\Bbb R)$ and $\varepsilon=0,1$.
It~remains to observe that $\varsigma_p^{-1}(v)\in S$ is equivalent to
$\big(\varepsilon-f(-p)\big)\langle
v,v\rangle-2fv+\varepsilon+f(-p)=0$.

\smallskip

{\bf2.11.}~Return to this exercise after studying the elements of
riemannian geometry. The vector $\langle-,v\rangle q\in\T_q\Bbb S^n$,
where $v\in q^\perp$, is tangent to the curve
$c(t):=q+tv\in V^\centerdot$ at $c(0)=q$. Since the definition of
$\varsigma_p$ in Exercise 2.8 works in some open neighbourhood of $q$
in $V^\centerdot$, we obtain
$\varsigma_p\langle-,v\rangle q=\displaystyle\frac{\big(1+\langle
q,p\rangle\big)v-\langle v,p\rangle(q+p)}{(1+\langle q,p\rangle)^2}$.
Consequently,
$\big\langle\varsigma_p\langle-,v_1\rangle
q,\varsigma_p\langle-,v_2\rangle
q\big\rangle=\displaystyle\frac{\langle v_1,v_2\rangle}{\big|1+\langle
q,p\rangle\big|^2}$
for $v_1,v_2\in q^\perp$.

\smallskip

{\bf3.3.2.}~Let $f\in C^1(U)$ and $p\in U$. By the mean value theorem,
for every sufficiently small $\varepsilon>0$, there exists
$\varepsilon'\in[0,\varepsilon]$ such that
$\displaystyle\frac{f(p+\varepsilon v)-f(p)}\varepsilon=
v_{p+\varepsilon' v}f$.
Hence,
$\displaystyle\frac{f(p+\varepsilon w+\varepsilon v)-f(p+\varepsilon
w)}\varepsilon=v_{p+\varepsilon w+\varepsilon'v}f$
for every sufficiently small $\varepsilon>0$ and a suitable
$\varepsilon'\in[0,\varepsilon]$. We obtain
$$\frac{f(p+\varepsilon v+\varepsilon
w)-f(p)}\varepsilon=v_{p+\varepsilon
w+\varepsilon'v}f+\frac{f(p+\varepsilon w)-f(p)}\varepsilon.$$
Since $v_qf$ is continuous in $q\in U$, it follows that
$(v+w)_pf=v_pf+w_pf$.

\smallskip

{\bf3.3.3.}~For some $p\in U\opensub M$ and $f\in C^\infty(U)$, we have
$g=f_p$. The map $\varphi:v\mapsto v_pf$ (this definition is correct
since it is independent of the choice of $f$ representing $g$; so, we
can write $v_pg$) is a $\Bbb K$-linear functional by Exercise 3.3.2. It
follows from the Leibniz rule that $v_ph=0$ for all $v\in V$ and
$h\in\frak m_p^2$, implying the uniqueness.

Let $b^i\in V$ be some linear basis in $V$ and let $\varphi_i\in V^*$
be the corresponding dual basis. Then, by the Newton-Leibniz formula,
$f(v)=f(p)+\sum_i\varphi_i(v-p)f_i(v-p)$, where
$f_i(w):=\int_0^1b_{p+tw}^if\,\text dt$ is a smooth function in $w$ for
$w$ sufficiently close to $0$. It remains to apply the same formulae to
the functions $f_i(v-p)$ in $v$.

\smallskip

{\bf3.3.4.}~Show first that the function $f:\Bbb R\to\Bbb R$ given by
the rule $f(x):=0$ for $x\le0$ and $f(x):=\exp(-\frac1x)$ for $x>0$
is smooth. Then, assuming that $V$ is Euclidean, note that
$g^{-1}(\Bbb R^+)$ is the open ball of radius $r$ centred at $c$, where
$g(v):=f\big(r^2-\langle v-c,v-c\rangle\big)$ for all $v\in V$.

\smallskip

{\bf3.6.6.}~There is a homomorphism $h:\Cal F_p\to(\Cal F|_S)_p$ given
by the rule $f_p\mapsto(f|_{S\cap U})_p$, where $f\in\Cal F(U)$ and
$p\in U\opensub M$. By the definition of induced structure, $h$ is
surjective. It remains to observe that
$\ker h=\big\{f_p\in\Cal F_p\mid f\in\Cal F(U),\ f(S\cap U)=0\text{
for some }p\in U\opensub M\big\}$.

\smallskip

{\bf4.3.3.}~Show that $\displaystyle\frac{\sin(2l_3)}{\sin\alpha_2}$ is
symmetric in $r_1,r_2,r_3$.

\smallskip

{\bf4.4.1.}~Use orthogonal coordinates.

\smallskip

{\bf4.5.2.}~If $\langle p_1,p_2\rangle\ne0$, we take
$W:=\Bbb Rp_1+\Bbb R\langle p_1,p_2\rangle p_2$. Let
$\langle p_1,p_2\rangle=0$ and $p_1\in W\le V$, where $\Bbb P_\Bbb KW$
is a geodesic. Then $\pi[p_1]p\in W$ and $p_2=\pi[p_1]p$ for a suitable
$p\in W$.

\smallskip

{\bf4.5.9.}~See Exercise 4.8.?.

\smallskip

{\bf4.8.2.}~Stealing something from the proof of Stolen Carlos' lemma
is useful but this does not suffice.

\smallskip

{\bf6.5.}~Using the induction on $\dim W$, decompose
$W=W'\oplus\Bbb K w$. Being ${W'}^\perp\cap(\Bbb Kw)^\perp$ the kernel
of the functional ${W'}^\perp\to\Bbb K$ given by the rule
$x\mapsto\langle x,w\rangle$, we have
$\dim W^\perp=\dim\big({W'}^\perp\cap(\Bbb
Kw)^\perp\big)\ge\dim{W'}^\perp-1$
by Exercise 6.3. The rest follows from $\dim W=\dim W'-1$ by induction.

\smallskip

{\bf6.6.}~$W\cap W^\perp$ is the kernel of the induced form on $W$.

\smallskip

{\bf6.7.}~Use $W\subset{W^\perp}^\perp$ and Exercise 6.6.

\smallskip

{\bf6.9.}~Assuming that $\langle v,v\rangle=0$ for all $v\in V$, we
obtain $\langle v_1+v_2,v_1+v_2\rangle=0$ and, hence,
$\Re\langle v_1,v_2\rangle=0$ for all $v_1,v_2\in V$. It remains to
apply the last identity to $iv_1,v_2$.

\smallskip

{\bf6.10.}~Using Exercises 6.8 and 6.9, we can find a nonisotropic
$w\in W^\perp$ and put $W':=W+\Bbb Kw$.

\smallskip

{\bf6.23.}~Consider $m=\dim(W\cap V^\perp)$ and apply Exercise 6.4.

\smallskip

{\bf6.24.}~Decompose $V$ into the orthogonal sum $V=V_-\oplus V_+$ of
subspaces of signatures $(n_-,0,0)$ and $(0,0,n_+)$. If, say,
$\dim W>n_-\le n_+$, then $W\cap V_+\ne0$. In order to construct a
subspace of dimension $\min(n_-,n_+)$ with the null induced form, use
the isotropic elements mentioned in Example 6.19.

\smallskip

{\bf6.25.}~Decompose $W$ into the orthogonal sum
$W=W_-\oplus W_0\oplus W_+$ of subspaces of signatures $(m_-,0,0)$,
$(0,m_0,0)$, and $(0,0,m_+)$. Decomposing
$V=(W_-+W_+)\oplus(W_-+W_+)^\perp$, notice that
$m_-\le n_-$, $m_+\le n_+$, and $W_0\le(W_-+W_+)^\perp$, where
$(W_-+W_+)^\perp$ has signature $(n_--m_-,0,n_+-m_+)$. Using Exercise
6.24, conclude that $m_0\le n_--m_-$ and $m_0\le n_+-m_+$. For
$m_-$, $m_0$, and $m_+$ that satisfy the above inequalities, construct
a subspace of signature $(m_-,m_0,m_+)$.

\enddocument